\documentclass[11pt,a4paper]{amsart}
\usepackage{amsthm, amsmath, amsfonts, amssymb,graphicx}
\usepackage{mathrsfs,url,color,framed}
\usepackage{caption}
\usepackage[utf8]{inputenc}
\usepackage[T1]{fontenc}

\newtheorem{defn}{Definition}[subsection]
\newtheorem{exa}[defn]{Example}
\newtheorem{nonexa}[defn]{Non-Example}
\newtheorem{exaterm}[defn]{Example and Terminology}
\newtheorem{thm}[defn]{Theorem}
\newtheorem{goal}[defn]{Goal}
\newtheorem{exe}[defn]{Exercise}
\newtheorem{prop}[defn]{Proposition}
\newtheorem{propdef}[defn]{Proposition-Definition}
\newtheorem{conv}[defn]{Convention}
\newtheorem{notn}[defn]{Notation}
\newtheorem{lem}[defn]{Lemma}
\newtheorem{cor}[defn]{Corollary}
\newtheorem{wish}[defn]{Wishes}
\newtheorem{rem}[defn]{Remark}
\newtheorem{add}[defn]{Addendum}

\newtheorem{tenthm}[defn]{Tentative Theorem}

\numberwithin{equation}{subsection}

\newcommand{\CC}{\mathbb{C}}

\newcommand{\NN}{\mathbb{N}}

\newcommand{\QQ}{\mathbb{Q}}
\newcommand{\RR}{\mathbb{R}}

\newcommand{\ZZ}{\mathbb{Z}}

\newcommand{\cB}{{\mathcal{B}}}

\newcommand{\cM}{{\mathcal{M}}}

\newcommand{\cO}{{\mathcal{O}}}

\newcommand{\atup}{\underline{a}}
\newcommand{\btup}{\underline{b}}
\newcommand{\ctup}{\underline{c}}
\newcommand{\xtup}{\underline{x}}
\newcommand{\Xtup}{\underline{X}}
\newcommand{\ytup}{\underline{y}}
\newcommand{\Ztup}{\underline{Z}}
\newcommand{\mot}{^{\mathrm{mot}}}
\newcommand{\motloc}{^{\mathrm{mot,loc}}}

\newcommand{\supp}{\operatorname{supp}}
\newcommand{\rv}{\operatorname{rv}}
\newcommand{\RV}{\operatorname{RV}}
\newcommand{\res}{\operatorname{res}}

\newcommand{\rtrdim}{\operatorname{rtrdim}}
\newcommand{\gr}{\operatorname{gr}}
\newcommand{\lin}{\operatorname{line}}

\newcommand{\Cbb}[1]{\CC(\!(#1)\!)}
\newcommand{\Rbb}[1]{\RR(\!(#1)\!)}
\newcommand{\Hahn}{\Rbb{t^{\QQ}}}
\newcommand{\Hahnk}{k(\!(t^{\Gamma})\!)}
\newcommand{\Rs}{(\RR^*)}
\newcommand{\Rsn}{\Rs^n}
\newcommand{\Tr}{\operatorname{Tr}}
\newcommand{\Gro}{K_0(\mathrm{Var}_{\CC})}

\newcommand{\mpp}[1]{\marginpar{\raggedright\scriptsize #1}}
\newcommand{\empp}[2]{\emph{#1}\mpp{#2}}
\newcommand{\emp}[1]{\empp{#1}{#1}}

\newcommand{\inclgr}[2]{%
\fbox{%
\includegraphics[page=#1,trim={#2},clip,scale=0.5]{pictures}
}%
}

\newcommand{\pict}[4][t]{%
\begin{figure}[#1]
\inclgr{#2}{#3}%
\caption{#4 \label{fig.#2}}
\end{figure}%
\def\figref{Figure~\ref{fig.#2}}%
}

\newbox\figleft

\newcommand{\pictl}[3]{%
\global\setbox\figleft=\hbox{\inclgr{#1}{#2}}%
\gdef\captleft{#3 \label{fig.#1}}%
\def\figref{Figure~\ref{fig.#1}}%
}

\newcommand{\pictr}[4][t]{%
\begin{figure}[#1]
\begin{minipage}{0.5\textwidth}
\centerline{\box\figleft}%
\caption{\captleft}
\end{minipage}%
\begin{minipage}{0.5\textwidth}
\centerline{\inclgr{#2}{#3}}%
\caption{#4 \label{fig.#2}}
\end{minipage}
\end{figure}%
\def\figref{Figure~\ref{fig.#2}}%
}

\definecolor{gruen}{RGB}{0,100,0}

\def\bg{\bgroup\color{gruen}}
\def\eg{\egroup}

\title{A Non-Archimedean Approach to Stratifications}

\author{Immanuel Halupczok}

\begin{document}

\captionsetup{margin=0pt}

\maketitle

In \cite{i.whit,iY.sts,iC.aas,iBW.can}, Bradley-Williams, Cubides Kovacsics, Yin and the author of the present notes developed a new approach to stratifications that uses (non-archimedean) valued fields as a tool to have a good notion of infinitesimal neighbourhoods of singularities. In the present lecture notes, I will explain this approach, mainly following the version obtained in collaboration with Bradley-Williams in \cite{iBW.can}. (If you prefer a much shorter introduction, I recommend the introduction of \cite{iBW.can}.)
I tried to keep the notes as accessible as possible. In various places, it is useful to have some rudimentary knowledge of model theory, as provided e.g.\ by the lecture notes of Cubides Kovacsics in the same volume, but I hope that even readers with no model theory knowledge at all will be able to understand most of the things.

Here is an overview: Our main goal is to come up with a new notion of stratifications that is more natural than the classical ones, and also stronger than the (very classical) Whitney stratifications. I will start by explaining what exactly we want (Section~\ref{s.goal}) and then explain in detail how to obtain such stratifications using valued fields (Sections~\ref{s.nsa} to \ref{s.back}); in particular, Section~\ref{s.nsa} is a very short crash course to (a baby version of) non-standard analysis. In Section~\ref{s.rt}, we will see that our stratifications are only the shadow of something much more powerful, which we call the riso-tree. I will then explain some ingredients of the main proofs (Section~\ref{s.def}), and in the last section, I present an application to Poincaré series.

These lecture notes were written during a course I gave in Bilbao in summer 2023.
I am grateful to Javier Fernandez de Bobadilla for organizing this wonderful summer school and for giving me the opportunity to give this course. I am also grateful to the audience, for showing so much interest in this topic and for asking intelligent and interesting questions, some of which made it into these lecture notes.


\section{The goal}
\label{s.goal}

\subsection{A wish list for stratifications}

Suppose we are given some kind of nice geometric set $X \subset \RR^n$. For example, let $X$ be given as a zero set of some polynomials:

\begin{defn}\label{defn.alg}
We call a subset $X \subset \RR^n$ \empp{algebraic}{algebraic set $X$} if it is of the form
\[
X = \{\atup \in \RR^n \mid f_1(\atup) = \dots = f_\ell(\atup) = 0\},
\]
for some polynomials $f_1, \dots, f_\ell \in \RR[x_1, \dots, x_n]$. In the following, we fix such an algebraic set $X$.
\end{defn}

Such a set is smooth almost everywhere. Our goal is to get some understanding of the singular points of $X$.

\pagebreak 

\begin{defn}
Call $X$ \emp{smooth} at a point $\atup \in X$ if there exists a neighbourhood $U \subset X$ of $\atup$ which is a $C^\infty$-manifold.
If $X$ it not smooth at $\atup$, we say that $X$ is \emp{singular} at $\atup$.
\end{defn}

\begin{exa}\label{exa.kein}
The sphere, defined by $x^2 + y^2 + z^2 - 1 = 0$, is smooth everywhere
\pictl{1}{5cm 7cm 3.5cm 2.5cm}{Example~\ref{exa.kein}}(\figref).
\end{exa}
\begin{exa}\label{exa.pt}
The double cone, defined by $x^2 - y^2 - z^2 = 0$ has a singularity at $0$
\pictr{2}{4cm 7cm 2.5cm 2.5cm}{Example~\ref{exa.pt}}(\figref).
\end{exa}

\begin{exa}\label{exa.gerade}
The set defined by $y^2 + z^3 - z^2 = 0$ has an entire singular line
\pictl{3}{4cm 7cm 3cm 2.5cm}{Example~\ref{exa.gerade}}(\figref).
\end{exa}

A question we want to ask is: Among the singular points, can there be some which are ``more singular'' than others?

\begin{exa}\label{exa.singsing}
The set of singular points of the set defined by $x(y^2 + z^3 - z^2) = 0$ itself has a singularity at the origin
\pictr[b]{4}{4cm 6cm 3cm 2.5cm}{Example~\ref{exa.singsing}}(\figref).
This one surely deserves to be called ``more singular''.
\end{exa}

But is the singularity of Example~\ref{exa.pt} more singular than the ones from 
Example~\ref{exa.gerade}? Here is another example:

\begin{exa}\label{exa.whit}
In the set $X$ defined by $y^2 + z^3 - x^2z^2 = 0$, the set of singular points is just a line, so if we take ``more singular'' as meaning ``singular point of the set of singular points'', then the origin is not more singular. However, the origin is special, namely: If we consider the fiber of $X$ at a fixed $x$-coordinate, then for any $x\ne 0$, the fiber has a little loop,
whereas for $x = 0$, it is a cusp curve
\pict{5}{4cm 8cm 4cm 4cm}{Example~\ref{exa.whit}}(\figref).
\end{exa}

My intuition of what we mean by ``more singular'' is the following:
Near the origin in Examples~\ref{exa.pt}, \ref{exa.singsing} and \ref{exa.whit}, the surface $X$ is not translation invariant at all, whereas near all the other singular points in 
Examples~\ref{exa.singsing}, \ref{exa.whit} and \ref{exa.gerade}, it is almost translation invariant in some direction. (And at smooth points, the surfaces are almost translation invariant even in two directions.) Our goal is to make this precise. In other words, we would like to define a partition $X$ into sets $S_0, \dots, S_{\dim X}$, where $S_d$ is the set of
points near which $X$ is almost translation invariant in $d$ directions. (Such a partition is called a \emp{stratification}.) Here are some properties one would like a stratification to have:

\begin{wish}\label{w}
For each $d$, we would like to have the following:
\begin{enumerate}
 \item $S_d$ is a smooth manifold\dots
 \item \dots of dimension $d$ (unless it is empty).
 \item $S_0 \cup \dots \cup S_d$ is topologically closed.
 \item $S_d$ can be defined by a finite Boolean combination of polynomial equations.
 \item This should ``capture as many anomalies as possible''. In particular, in Example~\ref{exa.whit}, the origin should be in $S_0$.
\end{enumerate}
\end{wish}

Some remarks about those wishes:
\begin{rem}\label{rem.wish}
\begin{itemize}
 \item If $X$ would be a fractal curve, then one would expect $S_0$ to be all of $X$ (i.e., not of dimension $0$ and not smooth). However, we do not expect this to happen if $X$ is an algebraic set.
 \item A reason for Wish~(3) is: if arbitrarily close to some point $\atup \in X$, there are some bad singularities, then $\atup$ itself should be considered as being at least as bad. (A maybe more important reason is: This is useful for applications.)
\item The reason for Wish~(4) is that we would like the sets $S_d$ to live in the same category as our original set $X$. However, we certainly cannot expect the $S_d$ to be algebraic in the sense of Definition~\ref{defn.alg}, since we might have to remove lower-dimensional subsets (like the origin from $S_1$ in  Example~\ref{exa.singsing}).
\end{itemize}
\end{rem}

Given that the $S_d$ are Boolean combinations of algebraic sets,
one could decide to allow $X$ to be a Boolean combination as well. Being a model theorist, I would like to go even one step further:


\begin{conv}\label{c.l}
For the remainder of these notes, we fix a first order language $L$ 
\mpp{the language $L$} and consider $\RR$ as a structure in that language. Different options for $L$ are possible; for simplicity, let us fix $L$ to be the language of ordered rings; and for convenience, in these notes, let us moreover add a constant for each element of $\RR$ to the language:\footnote{By adding those constants to $L$, any set defined by polynomials over $\RR$ becomes $L$-definable.} $L := L_{\mathrm{oring}} \cup \RR = \{0, 1, +, -, \cdot, <\} \cup \RR$.\mpp{the language $L_{\mathrm{oring}}$ of odered rings}
\end{conv}

Non-model theorists, do not despair:

\begin{rem}\label{r.tar}
By Tarski's quantifier elimination result, the $L$-definable subsets of $\RR^n$ are exactly the finite Boolean combinations of sets defined by polynomial equations $f = 0$ and inequations $f \ge 0$, for some polynomials $f \in \RR[x_1, \dots, x_n]$.
\end{rem}

With the language $L$ fixed, we can reformulate Wish~(4) in a more natural way as follows:

\begin{wish}\label{w2}
Assuming that $X$ is $L$-definable, we would like to have:
\begin{enumerate}
 \item[(4')] $S_d$ is $L$-definable.
\end{enumerate}
\end{wish}

(While one can certainly not expect to obtain nice stratifications for arbitrary languages $L$, everything does work for various other ``tame'' languages; see Subsection~\ref{s.gen}. Most importantly, everything also works over $\CC$, with $L = L_{\mathrm{ring}} \cup \CC = \{0, 1, +, -, \cdot\} \cup \CC$.)

One way to achieve Wishes (1)--(3) is as follows: If $X$ has dimension $m$, then
let $S_m$ be the set of all those points of $X$ which have a neighbourhood which is an $m$-dimensional smooth manifold. Then repeat the process with the remainder $X \setminus S_m$ to define $S_{m-1}, S_{m-2}, \dots, S_0$.
It is not difficult to verify that this definition also fulfills Wish (4), but we already saw that (5) fails (in Example~\ref{exa.whit}). Thus the main question is: How to make formal that the origin in Example~\ref{exa.whit} is special?

\subsection{Some (unsuccessful) trials to define almost translation invariance}

As explained after Example~\ref{exa.whit}, my point of view that the reason for the origin in that example to be special is that it has no neighbourhood $U$ such that $X \cap U$ is almost translation invariant. If \emp{$d$-trivial} would be some formal definition of ``almost translation invariant in direction of a $d$-dimensional space'', then we would define the sets $S_d$ as follows:


\begin{defn}\label{d.Sd}
Let $S_d$ be the set of those points $\atup \in X$ such that $X$ is $d$-trivial but not $(d+1)$-trivial near $\atup$.
\end{defn}

We could try to make the ``almost'' from ``almost translation invariant'' formal by saying that $X$ becomes translation invariant after applying a suitable kind of map $\alpha$:

\begin{defn}\label{defn.dtriv}
Given $\atup \in X$,
we say that $X$ is \emp{$\underline{\qquad}$-$d$-trivial} near $\atup$ if there exists a neighbourhood $U \subset \RR^n$ and a $\underline{\qquad}$-bijection
$\alpha\colon U \to U'\subset \RR^n$ such that $Y := \alpha(X \cap U)$ is $d$-dimensionally translation invariant within $U$\mpp{$d$-dimensionally translation invariant within $U$}, i.e., there exists a $d$-dimensional vector sub-space $V \subset \RR^n$ such that for any two points $\btup_1, \btup_2 \in U'$ whose difference lies in $V$, we have $\btup_1 \in Y \iff \btup_2 \in Y$
\pict{34}{3cm 10.5cm 5cm 3.5cm}{Definition~\ref{defn.dtriv}}(\figref).
\end{defn}

It remains to decide what to write into the gaps.
The weakest condition which seems to make sense is ``homeomorphic''.\mpp{homeomorphically $d$-trivial} While
I would certainly expect almost translation invariance to be a stronger condition than just translation invariance up to homeomorphism, at least this captures Wish~(5). However, Wish~(4) seems to fail quite badly: existence of a homeomorphism is a pretty non-algebraic condition.\footnote{Though to be honest, I do not know of a concrete example where it fails.}

At the other end of the strength scale, it seems natural to impose that $\alpha$ is a smooth map. \mpp{smoothly $d$-trivial}However, this fails completely; the following example is due to Whitney.

\begin{exa}\label{exa.crr}
Let $X \subset \RR^3$ be defined by $yz(y-z)(y-xz) = 0$
\pict[b]{6}{3.5cm 8cm 3.5cm 3.5cm}{Example~\ref{exa.crr}}(\figref).
Then $X$ is not smoothly $1$-trivial near any of the points $(a,0,0)$ of the $x$-axis, so that the entire $x$-axis would be contained in $S_0$ (contradicting Wish~(2)). Indeed, consider the derivative at $(a,0,0)$ of a smooth map
$\alpha\colon U \to U'$, for some neighbourhood $U'$ of $a$. We may without loss suppose that it sends the three planes $y = 0$, $z = 0$ and $(y-z) = 0$ to themselves; but then, by the cross ratio obstruction, it also sends the plane $y - az = 0$ to itself; thus,
the image of the surface $y-xz = 0$ under 
$\alpha$ is not translation invariant.
\end{exa}

This example shows that even imposing $\alpha$ to be $C^1$ is too strong.

It is possible to impose that $\alpha$ is bi-Lipschitz.\mpp{bi-Lipschitz $d$-trivial} Compared to homeomorphic $d$-triviality, this feels a bit closer, in strength, to what I would imagine as ``almost translation invariant'', but it also has the problem that Wish~(4) fails.

\subsection{The classical solution and our new approach}
\label{s.csni}

Before continuing our own quest for a good notion of stratification, let me mention the classical solution(s) to the problem: The idea to obtain Wish~(4) is to define the sets $S_d$ in a much more algebraic way (in terms of certain tangent spaces converging to certain other tangent spaces). One then proves \emph{a posteriori} that $X$ satisfies some version of $d$-triviality near points of $S_d$. Concretely:

\begin{enumerate}
 \item A very classical notion of stratification is the one by Whitney \cite{Whi.strat}.\mpp{Whitney stratifications} It has a reasonably simple definition in terms of tangent spaces (called \emp{Whitney's regularity conditions}), and it implies that $X$ is homeomorphically $d$-trivial near points of $S_d$.
 \item \mpp{Mostowski's Lipschitz stratifications}Mostowski's Lipschitz stratifications \cite{Mos.biLip} even yield bi-Lipschitz $d$-triviality near points of $S_d$, but they are extremely technical to define.
\end{enumerate}

Neither of these notions yields an if-and-only-if: There might be points $\atup \in S_d$ near which $X$ is $d'$-trivial for $d' > d$. Moreover, those stratifications are not even defined canonically: We only obtain the existence of a stratification of $X$, but which stratum some point $\atup \in X$ lies in might not directly tell something about that point.\footnote{Whitney stratifications can at least somehow be made canonical, but Mostowski's Lipschitz stratifications are inherently non-canonical; see the discussion following Example~\ref{e.trumpet}.}

\begin{goal}
Find a notion of $d$-triviality which can directly serve to define a stratification using Definition~\ref{d.Sd}. 
\end{goal}

The most difficult part of this quest is to get Wish~(4) to hold. In addition, we hope that our notion of $d$-triviality captures well our intuition of ``almost translation invariance''.

Thinking about that intuition more precisely,
maybe the approach from Definition~\ref{defn.dtriv} was doomed from the beginning. After all, this definition only says ``almost translation invariant on a small neighbourhood''; whereas what we should really express is the following:
\begin{itemize}
 \item[(*)]
The smaller the neighbourhood of $\atup$, the closer to translation invariant $X$ should be.
\end{itemize}

This sounds as if we would end up with a messy condition involving lots of epsilons and deltas. Instead of even trying to make this precise, let us avoid this by using
non-standard analysis. The idea is to replace our field $\RR$ by some suitable bigger field $\RR^*$ which contains infinitesimal elements. In $\RR^*$, one can then make sense of the following variant of (*):
\begin{itemize}
 \item[(**)] On the \emph{infinitesimal} neighbourhood of $\atup$ (consisting of all elements infinitesimally close to $\atup$), $X$ should be \emph{infinitesimally} close to translation invariant.
\end{itemize}

A key principle from non-standard analysis (I will simply call it the \emp{the nsa-principle}) says that statements about infinitesimals in $\RR^*$ are equivalent to limit statements in $\RR$. In particular, using this principle, (**) in $\RR^*$ should be equivalent to (*) in $\RR$.

Usually, the nsa-principle provides a recipe to translate definitions formulated using $\RR^*$ back to definitions using purely $\RR$, so once we formulated (**) precisely, we should also be able to also get a precise formulation of (*) (thereby avoiding the need for $\RR^*$). 
However, formally, the nsa-principle cannot be applied directly to the condition (**) we will obtain. Thus, I am unable to make (*) formal, and I have to stick to (**) (in $\RR^*$) as the definition of $d$-triviality. Nevertheless, one can draw many conclusions from $d$-triviality which illustrate that (**) really behaves like (*).
And the real surprise is that this definition will also make Wish~(4) come true.

Let me now make (**) precise. While I could simply state the definition right away,\footnote{Very impatient readers can just read Convention~\ref{c.R*} and Definitions~\ref{d.app}, \ref{d.ri} and \ref{d.rt}, and then go to Section~\ref{s.back}.} I will rather first spend some time introducing tools and ideas from non-standard analysis, to justify the definition and to illustrate how the nsa-principle applies.


\section{Non-standard analysis}
\label{s.nsa}

\subsection{Adding infinitesimal elements to $\RR$}

We want an ordered field extension $\RR^*$ of $\RR$ which contains non-zero \empp{infinitesimal elements}{infinitesimal elements}, i.e., elements $a \in \RR^* \setminus \{0\}$ satisfying $-r < a < r$ for every $r \in \RR_{>0}$. For the nsa-principle to hold, we have to be quite careful in the choice of $\RR^*$, but there is still a lot of choice. Let me first specify the most important condition on $\RR^*$. (For non-model theorists, a concrete such $\RR^*$ will be specified in the next subsection.)

\begin{conv}\label{conv.R*}
Consider $\RR$ as an $L$-structure (where $L$ is the language from Convention~\ref{c.l}) and fix $\RR^*$\mpp{an elementary extension $\RR^* \succ \RR$} to be an elementary extension of $\RR$ which contains (non-zero) infinitesimal elements.\footnote{Recall that ``elementary extension'' means that for every $L$-formula $\phi(\xtup)$ and every tuple $\atup \in \RR^n$, we have $\RR \models \phi(\atup) \iff \RR^* \models \phi(\atup)$. For $L \supset K$, the notation ``$L \succ K$'' means: $L$ is an elementary extension of $K$.}
\end{conv}

%

Such elementary extensions always exist (independently of the language $L$ we use). Usually, in non-standard analysis, one simply puts \emph{all} subsets of $\RR^n$ (for every $n$) as relations into the language, to make everything definable. However, we will, at some point, need some more precise control over $\RR^*$, so that we have to be much more careful in our choice of $L$.

\subsection{A concrete example of $\RR^*$: Hahn series}

Let me specify a concrete example of such an elementary extension $\RR^*$. Later, we will need $\RR^*$ to satisfy an additional (non-model theoretic) assumption, which is satisfied by the example I give, so let us right away fix $\RR^*$ like this.

\begin{conv}\label{c.R*}
Let $\RR^* := \Hahn$\mpp{a concrete $\RR^*$} be the field of Hahn series over $\RR$ with rational exponents. Recall that, for any field $k$ and any ordered abelian group $\Gamma$, the \empp{Hahn field}{the Hahn field $k(\!(t^\Gamma)\!)$} $k(\!(t^\Gamma)\!)$ is defined as follows:
\[
k(\!(t^\Gamma)\!) := \{a = \sum_{\lambda \in \Gamma} a_\lambda t^\lambda \mid a_\lambda \in k, \supp a\text{ is well-ordered} \},
\]
where the support of $a$ is defined as \mpp{the support $\supp a$}$\supp a := \{\lambda \in \Gamma \mid a_\lambda \ne 0\}$.
\end{conv}
The elements of $k(\!(t^\Gamma)\!)$ are formal power series, which are added and multiplied as the notation suggests. Multiplying them involves taking certain sums of products of coefficients; the condition about the support being well-ordered ensures that those sums are always finite.

To see that $\Hahn$ is an elementary extension of $\RR$, we use two ingredients (neither of which we will prove).

First ingredient:

\begin{prop}
$\Hahn$ is real closed\mpp{real closed fields} (meaning it is not algebraically closed, but it becomes algebraically closed after adjoining a square-root of $-1$).
\end{prop}

\begin{rem}
The field $\Hahn$ has a natural order, which can be defined in the following two (equivalent) ways.\mpp{$<$ on $\Hahn$}
\begin{itemize}
 \item $a \ge b$ if and only if $a - b$ is a square in $\Hahn$. (This definition works in any real closed field.)
 \item We extend the order on $\RR$ to $\Hahn$ by imposing that $t$ is infinitesimal and positive. In other words, $a = \sum_{\lambda \in \QQ} a_\lambda t^\lambda$ is positive if and only if its leading coefficient $a_{\min \supp a}$ is positive.
\end{itemize}
\end{rem}

Second ingredient:

\begin{thm}
If $K \supset \RR$ is real closed, then $K$ is an elementary extension of $\RR$.
\end{thm}
(This theorem follows easily from Tarski's quantifier elimination result mentioned in
Remark~\ref{r.tar}.)

\begin{rem}
The field of Puiseux series over $\RR$ is also real closed, so it could, for the moment, also serve as an elementary extension of $\RR$. However, as already hinted at above, we will later need our elementary extension to have an additional property. This will disqualify the field of Puiseux series.
\end{rem}

\subsection{The nsa-principle in action}
\label{s.nsaa}

Let us now see what an elementary extension $\RR^* \succ \RR$ is useful for. More precisely, let us prove some first instances of the nsa-principle mentioned in Subsection~\ref{s.csni}. Without further mentioning, we will extend various \mpp{basic notation on $\RR^*$}basic notations from $\RR$ to $\RR^*$ in a natural way. \mpp{$|\cdot|$ on $\RR^*$}For example, the absolute value $|a|$ of an element $a \in \RR^*$ is equal to $a$ if $a \ge 0$ and equal to $-a$ if $a < 0$. This is then used to define the maximum norm on $\Rsn$: \mpp{$\|\cdot\|_\infty$ on $\Rsn$}$\|(a_1, \dots, a_n)\|_\infty = \max_i |a_i|$.

Using that $\RR^*$ is an elementary extension of $\RR$,
one easily verifies that properties of those notations are inherited. For example, the triangle inequality for the maximum norm can be expressed as an $L$-formula, so since it holds in $\RR^n$, it also holds in $\Rsn$.

Suppose that $X \subset \RR^n$ is $L$-definable, say, $X = \{\atup \in \RR^n \mid \RR \models \phi(\atup)\}$, for some $L$-formula $\phi(\xtup)$. Let
\mpp{$X^*$}$X^* := \{\atup \in \Rsn \mid \RR^* \models \phi(\atup)\}$ be the set defined by the same formula in $\RR^*$. (For non-model theorists: By Remark~\ref{r.tar}, $X$ is defined using polynomial equations and inequations. Evaluate the same conditions in $\RR^*$ to get $X^*$.) As a first example of how the nsa-principle works, suppose that 
we want to express whether some point (let us just take the origin) lies in the topological closure $\bar X$ of $X$. Using $\RR^*$, this can be expressed as follows:

\begin{prop}\label{p.cl}
$0 \in \bar X \iff X^*$ contains infinitesimal elements.
\end{prop}

\begin{proof}
$\Rightarrow$: We have $\RR \models \forall \epsilon > 0\colon \exists \atup \in X \colon \|\atup \|_\infty < \epsilon$. Since $\RR^*$ is an elementary extension, we thus also have $\RR^* \models \forall \epsilon > 0\colon \exists \atup \in X^* \colon \|\atup \|_\infty < \epsilon$. Now simply choose $\epsilon$ to be infinitesimal.

$\Leftarrow$: Let $\epsilon \in \RR_{>0}$ be given; we need to show that $\RR \models \exists \atup \in X \colon \|\atup \|_\infty < \epsilon$. Using $\RR^* \succ \RR$ again, this is equivalent to $\RR^* \models \exists \atup \in X^* \colon \|\atup \|_\infty < \epsilon$. This last statement is true by assumption: The norm of our infinitesimal $\atup \in X^*$ is in particular smaller than the (non-infinitesimal) $\epsilon$.
\end{proof}

\begin{exe}
Formulate a similar criterion for the origin lying in the topological interior of $X$.
\end{exe}

%
%

Next, suppose that $f\colon \RR \to \RR$ is an $L$-definable function. Then the formula which defines the graph of $f$ also defines the graph of a function from $\RR^*$ to $\RR^*$ (exercise); we denote that function by $f^*$. Intuitively, continuity means: if $x$ changes only a tiny bit, then also $f(x)$ only changes a tiny bit. In $\RR^*$, this intuition becomes formal:

\begin{prop}\label{p.cont}
The function $f$ is continuous at $0$ if and only if for every infinitesimal $x \in \RR^*$, $f^*(x) - f^*(0)$ is infinitesimal.
\end{prop}

The proof method is similar to that of Proposition~\ref{p.cl}, just slightly more involved because we now have more quantifier alterations. I nevertheless leave it as an exercise.

As a last example of the nsa-principle, let us express derivatives using $\RR^*$. I want to consider a slightly stronger notion of derivative than usual, namely the \emph{strict} derivative, since the corresponding statement in $\RR^*$ will be closely related to the notion of $d$-triviality we will introduce.

\begin{defn}
A function $f\colon \RR \to \RR$ is said to be \empp{strictly differentiable}{strict differentiability} at some point $a \in \RR$, with (strict) derivative $c \in \RR$, if
\[
\lim_{\substack{x_1,x_2 \to a\\x_1 \ne x_2}}\frac{f(x_1) - f(x_2)}{x_1 - x_2} = c.
\]
\end{defn}
(The difference to the usual notion of derivative is that we let both points vary instead of fixing $x_2$ to be $a$.)


By mimicking Proposition~\ref{p.cont}, the reader should be able to guess that we have the following equivalence (which, by the way, can easily be deduced from a multi-variable version of Proposition~\ref{p.cont}):

\begin{prop}\label{prop.std}
An $L$-definable function $f\colon \RR \to \RR$ has strict derivative $c$ at $0$ if and only if, for every pair of infinitesimals $x_1, x_2 \in \RR^*$ with $x_1 \ne x_2$,
\begin{equation}\label{eq.std}
\frac{f^*(x_1) - f^*(x_2)}{x_1 - x_2} - c \text{ is infinitesimal.} 
\end{equation}
\end{prop}

Condition~\eqref{eq.std} can be considered as saying that the graph of $f^*$ restricted to the infinitesimals is almost a straight line
\pict{7}{1.5cm 10cm 6cm 4cm}{Condition~\eqref{eq.std}}(\figref, left hand side). Note that if we would use the corresponding condition for normal derivatives instead, it would be an almost straight line only in a weaker sense (\figref, right hand side).


\section{$\RR^*$ as a valued field}

\subsection{Defining a valuation on $\RR^*$}

It will be handy to have some more tools to deal with the notion of infinitesimality. This is achieved by noting that $\RR^*$ comes with a natural (Krull) valuation. (If we choose $\RR^* = \Hahn$, we already have a valuation. Let me nevertheless define it in general.)

\begin{notn}
Let $\cO$\mpp{the valuation ring $\cO$} be the set of \empp{finite elements}{finite element} of $\RR^*$, i.e., elements $a \in \RR^*$ for which there exists an $r \in \RR_{>0}$ with $-r < a < r$.
Let $\cM$\mpp{the maximal ideal $\cM$} be the set of infinitesimal elements of $\RR^*$.
\end{notn}

We leave it to the reader to verify that $\cO$ is a valuation sub-ring of $\RR^*$ and that $\cM$ is its maximal ideal. We then define the corresponding valuation as usual:

\begin{notn}
Let \mpp{the value group $\Gamma$}$\Gamma := (\RR^*)^\times/\cO^\times$ be the quotient group, but written additively. Let \mpp{the valuation $v$}$v\colon  (\RR^*)^\times \to \Gamma$ be the canonical map, extended by $v(0) := \infty$. Put the order on $\Gamma$ defined by $v(a) \ge v(b)$ if $\frac ab \in \cO$.
\end{notn}

One easily verifies that the residue field $\cO/\cM$ is naturally isomorphic to $\RR$ (exercise).

\begin{notn}
We denote the residue map by \mpp{the residue map $\res$}$\res\colon \cO \to \RR$. (It sends $a \in \cO$ to the unique real number number infinitesimally close to $a$.)
\end{notn}

Intuitively, $v(a)$ is the negative of the order of magnitude\footnote{Looking up some latin suggests the term ``order of parvitude'' for this negative.} of $a$: We have $v(a) > v(b)$ if and only if $|a| < r\cdot |b|$ for every $r \in \RR_{>0}$.

\begin{exa}
In the case $\RR^* = \Hahn$, we recover the usual valuation:
$\Gamma$ can be identified with $\QQ$, and $v(a) = \min \supp a$.
\end{exa}

\subsection{Valuations of tuples and valuative balls}

Let me now fix some more notation. Let us do this in more generality:

\begin{conv}\label{c.K}
Let $K$ be an arbitrary Krull valued field, with valuation ring $\cO$, maximal ideal $\cM$, valuation $v\colon K \to \Gamma \cup \{\infty\}$ and residue field $k$. 
\end{conv}

Firstly, we define the valuation of a tuple. Intuitively, a tuple has small order of magnitude if all its entries have; this justifies the following definition.

\begin{defn}
\mpp{$v((a_1, \dots, a_n))$}
We define the valuation of a tuple $\atup = (a_1, \dots, a_n) \in K^n$ to be $v(\atup) := \min \{v(a_1), \dots, v(a_n)\}$.
\end{defn}

Note that in the case $K = \RR^*$, the valuation $v(\atup)$ is
not only equal to $v(\|\atup\|_\infty)$, but also to
the valuation of any other norm of $\atup$, e.g., to $v(\|\atup\|_2)$ (exercise).

Next, we introduce the notion of valuative balls:
\begin{defn}
Given $\atup \in K^n$ and $\lambda \in \Gamma$, we define a \empp{closed valuative ball}{closed valuative balls $B(\atup, \ge \lambda)$} to be a set of the form $B(\atup, \ge \lambda) := \{\btup \in K^n \mid v(\btup - \atup) \ge \lambda\}$. An \empp{open valuative ball}{open valuative balls $B(\atup, >\lambda)$} is a set of the form $B(\atup, > \lambda) := \{\btup \in K^n \mid v(\btup - \atup) > \lambda\}$. In both cases, we call $\lambda$ the \empp{(valuative) radius}{the (valuative) radius} of the ball.
\end{defn}

For the intuition, it is important to keep in mind that in $\Rsn$, there are much less valuative balls than metric balls. For example, $\cM^2$ is a very small ball (only infinitesimals), whereas the next bigger valuative ball is already $\cO^2$. In contrast, there are many metric balls between $\cM^2$ and $\cO^2$
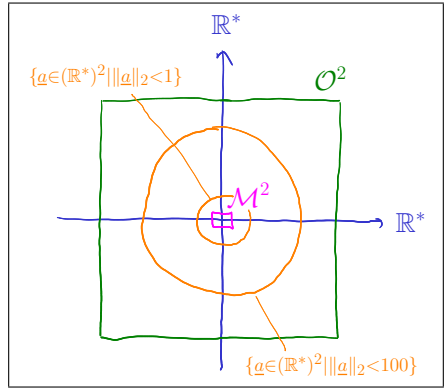
\pict{8}{4cm 7cm 4cm 2cm}{valuative and metric balls}(\figref).

\begin{exa}\label{exa.ball}
If $K = \Hahnk$, then we can write a closed valuative ball $B \subset K^n$ in the form $B = \sum_{\mu < \lambda} \atup_\mu t^\mu + t^\lambda\cdot \cO^n$ (where $t^\lambda\cdot \cO^n$ consists of the series having only $t$-exponents at least $\lambda$).
\end{exa}

\subsection{Hahn fields are spherically complete}

Let us take the opportunity to note that every Hahn field $\Hahnk$ is complete in the following strong sense. (This will later be needed.)

\begin{prop}\label{p.hsp}
$\Hahnk$ is spherically complete, i.e., every nested chain of valuative balls $(B_i)_{i \in I}$ has non-empty intersection: $\bigcap_{i \in I} B_i \ne \emptyset$.
(Nested chain means: $B_i \subset B_j$ or $B_j \subset B_i$ for every $i,j \in I$.)
\end{prop}
(For the notion of spherical completeness, it is irrelevant if we require the $B_i$ to be open or closed, or if we allow both.)

A somewhat counter-intuitive phenomenon is that even if the valuative radius of the $B_i$ does not go to $\infty$, their intersection could be empty:

\begin{nonexa}
The field $K := \bigcup_{n \ge 1}\Rbb{t^{1/n}}$ of Puiseux series\mpp{the field of Puiseux series} is not spherically complete:
Set
\[
B_0 := t^{1/2}\cO,\quad B_1 := t^{1/2}+t^{3/4}\cO,\quad B_2 := t^{1/2}+t^{3/4} + t^{7/8}\cO, \quad\dots.
\]
Clearly, $B_0 \supset B_1 \supset B_2 \supset \dots$, but a power series in the intersection $\bigcap_{i\in \NN} B_i$ would have to start with $b := t^{1/2} + t^{3/4} + t^{7/8}+\dots$. This cannot be an element of $K$, since the denominators in the $t$-exponents of $b$ are not bounded. (On the other hand, $b$ is an element of $\Hahn$.)
\end{nonexa}

One way to prove Proposition~\ref{p.hsp} goes as follows: Write each ball as in Example~\ref{exa.ball}. The balls being nested means that the partial series $\sum_{\mu < \lambda} \atup_\mu t^\mu$ agree on the terms they have in common. To find an element in the intersection, just take the series consisting of all the terms appearing in all of those partial series. (I leave the details to the reader.)

\subsection{Being almost equal -- the leading term structures $\RV^{(n)}$}

It will be useful to have a notion of ``being almost equal'', which means: equal up to something of smaller order of magnitude.

\begin{defn}\label{d.app}
Given $\atup, \btup \in K^n$, set\mpp{$\approx$}
\[\atup \approx \btup
:\iff v(\atup - \btup) > v(\atup) \quad\vee \quad \atup = \btup = 0.\]
We define \mpp{$\RV^{(n)}$}$\RV^{(n)} := K^n/\mathord{\approx}$ to be the quotient by that equivalence relation, and we denote the natural map to the quotient by \mpp{$\rv$}$\rv\colon K^n \to \RV^{(n)}$. Instead of $\RV^{(1)}$, one usually writes \mpp{$\RV$}$\RV$.
\end{defn}
\pict{9}{5cm 7cm 4cm 2cm}{Definition~\ref{d.app}}
In \figref, each of the red squares is one $\approx$-equivalence relation: Those in $\cO^2 \setminus \cM^2$ are just infinitesimal neighbourhoods of points; those within $\cM^2$ are much smaller. (And equivalence classes outside of $\cO^2$ would be too big to draw.)

\begin{exe}
Verify that $v(\atup - \btup) > v(\atup)$ implies $v(\atup) = v(\btup)$ and use this to deduce that the above relation $\approx$ is indeed symmetric.
\end{exe}

\begin{exaterm}
Suppose that $K = \Hahnk$, and
consider $\atup = \sum_\mu \atup_\mu t^\mu$ and $\btup = \sum_\mu \btup_\mu t^\mu$ in $K^n$. Then $\atup \approx \btup$ implies $v(\atup) = v(\btup) =: \lambda$, so the two series have leading terms with the same $t$-exponent $\lambda$. The condition $v(\atup - \btup) > \lambda$ is equivalent to additionally imposing $\atup_\lambda = \btup_\lambda$. In other words, $\atup \approx \btup$ exactly says that $\atup$ and $\btup$ have the same leading term. Therefore, $\RV$ is called the \empp{leading term structure}{the leading term structure}. (In the notation $\RV$, the V stands for valuation, and the R stands for ``residue'' and refers to the leading coefficient, which is an element of the residue field.)
\end{exaterm}

\begin{rem}
Note that $\RV^{(n)}$ is not the same as the Cartesian power $\RV^n$. Indeed consider $\atup = (1, 0), \btup = (1, t) \in \Hahn$. We have 
$\rv(\atup) = \rv(\btup)$ in $\RV^{(2)}$,
but $(\rv(1), \rv(0)) \ne (\rv(1), \rv(t))$ in $\RV^2$. (However, it is not difficult to verify that there is a canonical map $\RV^n \to \RV^{(n)}$.)
\end{rem}

To get some practice using $\approx$ and $\rv$, let us reformulate Condition~\eqref{eq.std} from Proposition~\ref{prop.std} (which is equivalent to an $L$-definable function $f\colon \RR \to \RR$ having strict derivative $c$ at $0$): Given $x_1, x_2 \in \cM$, set $y_i := f^*(x_i)$. We claim that the condition is equivalent to
\begin{equation}\label{eq.std2}
(x_1, y_1) - (x_2, y_2) \approx (x_1, cx_1) - (x_2, cx_2),
\end{equation}
which is another way to say that the graph of $f$ looks very much like the graph of $x \mapsto cx$
\pict{10}{1.5cm 9cm 6cm 4cm}{Condition~\eqref{eq.std2}}(\figref).

Indeed the right hand side of this new condition has
valuation
\[
\max\{v(x_1 - x_2), v(c(x_1 - x_2))\} = v(x_1 - x_2)
\]
(since $v(c) \ge 0$), and in the difference of both sides of \eqref{eq.std2}, the $x$-coordinates cancel, so that the valuation of that difference is
\[
v(y_1 - y_2 - (cx_1 - cx_2)).
\]
Thus \eqref{eq.std2} is (by definition of $\approx$) equivalent to
\begin{equation}\label{eq.std3}
 v(y_1 - y_2 - (cx_1 - cx_2)) > v(x_1 - x_2),
\end{equation}
which, after dividing both terms inside the $v(\cdot)$ by $x_1 - x_2$, yields exactly \eqref{eq.std}.

\section{Risometries}

\subsection{Maps which almost preserve translation invariance}

Recall that to define our stratifications, we want to come up with a notion of ``infinitesimally close to translation invariant''. We will do this using the idea from Definition~\ref{defn.dtriv}, namely by imposing that our set becomes translation invariant after applying a certain kind of map $\alpha$. The maps $\alpha$ we want to allow are ``risometries'', which we will define now.
The definition makes sense in any valued field $K$, so let me formulate it in this general context (see Convention~\ref{c.K}).
However, to prove properties of that notion, we will soon need to impose that $K$ is spherically complete, and maybe also that it is of equi-characteristic $0$, so to be on the safe side, let us impose those conditions.

Later, we will apply this when $K$ is an elementary extension of $\RR$. At that point, we will need a spherically complete elementary extension, so we will set $K := \RR^* := \Hahn$ (as in Convention~\ref{c.R*}).

\begin{defn}\label{d.ri}
A bijection $\alpha\colon U \to U'$ between subsets $U, U' \subset K^n$ is called a \empp{risometry}{risometries} if, for every $\atup_1, \atup_2 \in U$, we have:
\[
\alpha(\atup_1) - \alpha(\atup_2) \approx \atup_1 - \atup_2
\]
\pict{11}{3cm 8cm 4cm 5cm}{Definition~\ref{d.ri}}(\figref).
\end{defn}

\begin{rem}\label{rem.iso}
A risometry is in particular a valuative isometry, i.e., it satisfies 
$v(\alpha(a_1) - \alpha(a_2)) = v(a_1 - a_2)$.
In particular, risometries are continuous.
\end{rem}

The term ``risometry'' comes from the fact that its definition is obtained by taking the definition of isometry and inserting some ``r'': $\rv(\alpha(a_1) - \alpha(a_2)) = \rv(a_1 - a_2)$.

\begin{rem}
It is easy to see that the set of risometries from a set $U$ to itself forms a group under composition.
\end{rem}

\subsection{Why we need spherical completeness}

We will almost only be interested in the case where the domain $U$ of the risometry in Definition~\ref{d.ri} is a valuative ball, e.g.\ $U = B(\atup, >\lambda)$ for some $\atup \in K^n$ and some $\lambda \in \Gamma$. (The following discussion works equally for closed valuative balls.)
%
Using Remark~\ref{rem.iso}, one easily sees that a risometry $\alpha$ sends such a $U$ to a subset of $B(\alpha(\atup), >\lambda)$. At some point, it will be important that the image $\alpha(U)$ is always equal to the entire ball $B(\alpha(\atup), >\lambda)$ (see Exercise~\ref{exe.c1t} (1)). This is where we need the assumption that $K$ is spherically complete:

\begin{prop}\label{prop.surj}
Suppose that $K$ is a spherically complete valued field. Then for any risometry
$\alpha\colon U \to U' \subset K^n$ with $U = B(\atup, >\lambda)$, we have $U' = B(\alpha(\atup), >\lambda)$.
\end{prop}

\begin{proof}[Sketch of proof of Proposition~\ref{prop.surj}]
We already know that $U' \subset B(\alpha(\atup), >\lambda)$, so we just need to prove surjectivity of $\alpha$ onto $B(\alpha(\atup), >\lambda)$. Thus, let $c \in B(\alpha(\atup), >\lambda)$ be given and let us try to find a preimage.

We will do so by approximation. To this end, let us say that $b \in U$ has ``quality $\lambda$'' if $v(\alpha(b) - c) = \lambda$. We start by showing that ``quality can always be improved'': Given $b$ of quality $\lambda$, we find a $b'$ of quality $\lambda' > \lambda$ as follows
\pict[b]{12}{3cm 8cm 4cm 5cm}{Proof of Proposition~\ref{prop.surj}}(\figref):

Set $b' := b + (c - \alpha(b))$. Then, since $\alpha$ is a risometry,
we have $\alpha(b') - \alpha(b) \approx b' - b = c - \alpha(b)$,
which means that $v((\alpha(b') - \alpha(b)) - (c - \alpha(b))) > v(c - \alpha(b)) = \lambda$. On the left hand side, $\alpha(b)$ cancels, and we obtain the quality of $b'$.

If the value group would be $\ZZ$, then by repeatedly improving the quality as above, we would obtain a sequence $(b_i)_{i \in \NN}$ whose quality goes to $\infty$, so that the limit would be a preimage of $c$. (Here, we use completeness to have that the limit exists and continuity of $\alpha$ to deduce that the limit $\lim_i b_i$ is sent to $\lim_i c_i$.) For other value groups $\Gamma$, we can still make the argument work, using spherical completeness. Here is a sketch:

Let $S$ be the set of those $\lambda \in \Gamma$
which do arise as qualities of elements of $U$, and for each $\lambda \in S$, fix a $b_\lambda \in U$ of that quality.

Case 1: $S$ has a maximum $\mu$. Then applying the quality-can-always-be-improved argument to $b_\mu$ yields a contradiction to $\mu$ being maximal.

Case 2: $S$ has no maximum. Set $B_\lambda := B(b_\lambda, \ge \lambda)$ for $\lambda \in S$. I leave it as an exercise to verify that $B_\lambda$ is exactly the set of elements of quality at least $\lambda$. (This only uses that $\alpha$ is a valuative isometry.) From this, one deduces that the balls $B_\lambda$ form a nested chain, so by spherical completeness, they have a non-empty intersection $B$. Any element $b \in B$ has quality bigger than every $\lambda \in S$. The only way for this to not be a contradiction is that $b$ has quality $\infty$ (and $B = \{b\}$), meaning that $\alpha(b) = c$ and that we are done with the proof.
\end{proof}

\subsection{Exercises about and examples of risometries}

To get a better intuition about risometries, here is a little warm-up exercise:
\begin{exe}\label{e.rires}
Suppose that $\alpha \colon \cO^n \to \cO^n$ is a risometry.
\begin{enumerate}
 \item Show that
 $\alpha$ induces a map $\bar\alpha\colon k^n \to k^n$, where $k$ denotes the residue field of $K$. (I.e., we want that $\res(\alpha(\atup)) = \bar\alpha(\res(\atup))$ for all $\atup \in\cO^n$
 \pict[b]{13}{3cm 5.5cm 4cm 5cm}{Exercise~\ref{e.rires}}(\figref).
 \item Show that this induced map is just a translation, i.e., $\alpha(\bar\atup) = \bar\atup + \bar\btup$ for some fixed $\bar\btup \in k^n$.
 \item Is this an if-and-only-if condition? I.e.: If  $\alpha \colon \cO^n \to \cO^n$ is a map inducing a translation at the level of the residue field, does this imply that $\alpha$ is already a risometry?
 \item Bonus question:\footnote{Not useful for the present notes, but a fun exercise to get used to elementary extensions.} How does the answer to (3) change if $K = \RR^*$ and we require $\alpha$ to be definable in the language $L$ from Convention~\ref{c.l}?
\end{enumerate}
\end{exe}

We will be interested in whether a risometry can make a set translation invariant. In the next exercise, let us analyse precisely what this means in an easy special case.

\begin{exe}\label{exe.c1t}
\begin{enumerate}
 \item Suppose that $\alpha \colon \cM^2 \to \cM^2$ is a risometry. Show that the preimage $Z := \alpha^{-1}(\cM \times \{0\})$ of the horizontal line is the graph of a function (which we shall denote by $g$).
 
 Hint 1: First show that each fiber $Z_x := \{y \in \cM \mid (x,y) \in Z\}$ consists of at most one point. (For this, you only need to use the risometry condition.)
 
 Hint 2: To show that each fiber is non-empty, it may be useful to consider the map $\beta\colon \cM \to \cM$ sending $x$ to $\pi(\alpha^{-1}(x,0))$. Can you show that it is a risometry? And can you then see the real reason we need Proposition~\ref{prop.surj}?
 \item Show that this function $g$ from (1) satisfies the following condition:
\begin{equation}\label{eq.c1t}
\forall x_1, x_2 \in \cM\colon (x_1 \ne x_2 \Rightarrow v(g(x_1) - g(x_2)) > v(x_1 - x_2)).
\end{equation}
 \item Now let us do the converse: Show that if $g\colon \cM \to \cM$ is a function satisfying \eqref{eq.c1t}, then there exists a risometry $\alpha'\colon \cM^2 \to \cM^2$ sending the graph of $g$ to $\cM \times \{0\}$.
 
 Hint: First define $\alpha'$ on the graph of $g$ itself in the most natural possible way. Then extend the definition of $\alpha'$ to all of $\cM^2$, again in the most natural possible way.
\end{enumerate}
\end{exe}

\emph{A priori}, asking whether there exists a risometry sending a set $Z \subset \cM^2$ to $\cM \times \{0\}$ is a condition which model theorists do not like, since one has an existential quantifier over all possible risometries. However, in the exercise, we just saw that this is equivalent to a condition with quantifiers running only over (subsets of Cartesian powers of) the valued field:
\begin{equation}\label{eq.tf}
\begin{aligned}
&\forall x \in \cM\colon \exists^{=1} y \in \cM\colon (x,y) \in Z \,\,\wedge\\
&\forall (x_1,y_1), (x_2, y_2) \in Z: (x_1 \ne x_2 \Rightarrow v(y_1 - y_2) > v(x_1 - x_2)) 
\end{aligned}
\end{equation}
Results like that will be crucial for Wish~\ref{w2}~(4') concerning definability of the sets $S_d$.

Let us look at another consequence of Exercise~\ref{exe.c1t}. Firstly, note that it generalizes as follows:
Given any $c \in \cO$, a subset $Z \subset \cM^2$ can be sent to $\{(x, cx) \mid x \in \cM\}$ by a risometry if and only $Z$ is the graph of a function $g\colon \cM \to \cM$ which satisfies, for any distinct $x_1, x_2 \in \cM$: 
\begin{equation}\label{eq.c1t2}
 v(g(x_1) - g(x_2) - c\cdot(x_1 - x_2)) > v(x_1 - x_2)).
\end{equation}
This condition is exactly \eqref{eq.std3}, so we deduce the following connection to strict differentiability: Fix an $L$-definable set $X \subset \RR^2$, and suppose for the moment that $X$ contains the origin.

\begin{lem}\label{l.c1}
If there exists a (classical) neighbourhood $U \subset \RR^2$ of $(0,0)$ such that $X \cap U$ is the graph of a function which is strictly differentiable at $0$, then there exists a risometry $\alpha \colon \cM^2 \to \cM^2$ such that the image $\alpha(X^* \cap \cM^2)$ is a straight line
\pict{14}{5cm 9cm 5.5cm 4cm}{Lemma~\ref{l.c1}}(\figref).
\end{lem}

We claim that the converse of the lemma also holds -- except that we might have to swap coordinates if the image under the risometry is (almost) vertical. Here is the precise formulation of the converse:

\begin{lem}
If there exists a risometry $\alpha \colon \cM^2 \to \cM^2$ such that the image $\alpha(X^* \cap \cM^2)$ is a straight line $Z'$, then there exists a (classical) neighbourhood $U \subset \RR^2$ of $(0,0)$ such that, after possibly swapping coordinates, $X \cap U$ is the graph of a function which is strictly differentiable at $0$.
\end{lem}

\begin{proof}[Sketch of proof]
We may assume that the slope $c$ of $Z'$ has valuation $v(c) \ge 0$; if not, swap coordinates. Thus $X^* \cap \cM^2$ is the graph of a function $g$ satisfying \eqref{eq.c1t2}, which is almost what we want; the only problem left is that to conclude using our observations about strict differentiability, we would need $g$ to be equal to $f^*$, for some $L$-definable function $f$ on some small (classical) interval $(-\epsilon,\epsilon) \subset \RR$. This can be obtained using similar non-standard-analysis methods as in Subsection~\ref{s.nsaa}: In $\RR^*$, it is true that $X^* \cap U$ is the graph of a function, for $U$ sufficiently small; therefore, this also holds in $\RR$. (However, here, I swept some of subtleties under the carpet: It needs a bit of fiddling to find the right neighbourhoods $U$ to work with.)

Once we have a $U \subset \RR^2$ such that $X \cap U$ is the graph of $f$, we can finish as planned: Since $f^*$ satisfies \eqref{eq.c1t2}, $f$ is strictly differentiable at $0$.
\end{proof}

Applying the two previous lemmas to every point of $X$, we obtain the following result, which illustrates quite well how risometries in $\RR^*$ are useful in $\RR$.

\begin{prop}\label{p.mnf1}
An $L$-definable subset $X \subset \RR^2$ is a $1$-dimensional (strictly) $C^1$-sub-manifold of $\RR^2$ if and only if, for every $\atup \in X$, there exists a risometry $\alpha \colon B(\atup, >0) \to B(\atup, >0)$ such that $\alpha(X^* \cap B(\atup, >0))$ is a straight line.
\end{prop}

In the proposition, we can drop the word ``strictly'': if a function is $C^1$ on an open set $U \subset \RR^n$, then it is automatically strictly $C^1$ on $U$.

And as the reader will have guessed, similar statements also hold in higher dimension.

\section{Riso-triviality}

\subsection{Almost-translation-invariance: the definition}

Now we are ready to formulate our notion of ``infinitesimally close to translation invariant''. According to the terminology of Definition~\ref{defn.dtriv}, it would be called ``risometrically $d$-trivial'', but I will just write ``$d$-riso-trivial''.

We continue working in an arbitrary valued field $K$ of equi-characteristic $0$ which is spherically complete, e.g.\ $K = \Hahnk$.\footnote{Actually, every spherically complete valued field of equi-characteristic $0$ is isomorphic to such a Hahn field \cite{Kap.maxValFlds}, so instead of ``e.g.'', it might be more appropriate to write ``i.e.''.}

At some point, it will be useful to have the notion of $d$-riso-triviality not just for a single subset of $K^n$, but for a tuple of sets, so I right away formulate the definition for such tuples. Moreover, I use the opportunity to introduce some more terminology.

\begin{defn}\label{d.rt}
Let $\Ztup := (Z_1, \dots, Z_\ell)$ be a tuple of arbitrary sets $Z_i \subset K^n$ and let
$B \subset K^n$ be a valuative ball.
\begin{enumerate}
 \item Given a sub-vector space $V \subset K^n$, we call $\Ztup$ \emp{$V$-riso-trivial} on $B$ if there exists a risometry $\alpha\colon B \to B' \subset K^n$ such that each of the sets $\alpha(Z_1 \cap B), \dots, \alpha(Z_\ell \cap B)$ is translation invariant in the direction of $V$ (within $B'$).
\item We call $\Ztup$ \emp{$d$-riso-trivial} on $B$ (for $0 \le d \le n$) if there exists a $d$-dimensional $V \subset K^n$ such that $\Ztup$ is $V$-riso-trivial on $B$.
\item The \empp{riso-triviality dimension}{the riso-triviality dimension $\rtrdim_B Z$} $\rtrdim_B(\Ztup)$ of $\Ztup$ on $B$ is the maximal $d$ such that $\Ztup$ is \emph{$d$-riso-trivial} on $B$.\footnote{In \cite{iBW.can}, this dimension is denoted by $\dim \operatorname{rtsp}_B(\Ztup)$. Since I do not want to introduce $\operatorname{rtsp}_B(\Ztup)$ in these notes, I use the \emph{ad hoc} notation $\rtrdim_B(\Ztup)$. If you became curious about $\operatorname{rtsp}_B(\Ztup)$, you can solve Exercise~\ref{e.rtsp}\label{f.rtsp}.}
\end{enumerate}
\end{defn}

Let us start considering some trivial examples:
\begin{itemize}
 \item 
Obviously, $0$-riso-triviality always holds.
\item
For a set $Z \subset K^n$ to be $n$-riso-trivial on a ball $B \subset K^n$, we need $\alpha(Z \cap B)$ to be either empty or all of $B$, and hence also $Z \cap B$ itself has to be either empty or all of $B$.
\end{itemize}

A first non-trivial example has already been studied in Exercise~\ref{exe.c1t}: A graph of a function $g\colon \cM \to \cM$ is $(K \times \{0\})$-riso-trivial if and only if it satisfies \eqref{eq.c1t}.

Do we really need to introduce a notion for tuples of sets?
Clearly, for $(Z_1,Z_2)$ to be $d$-riso-trivial, it does not suffice that both $Z_1$ and $Z_2$ are $d$-riso-trivial. Note however that this is not even true if $Z_1$ and $Z_2$ are both $V$-riso-trivial
for the same sub-vector space $V \subset K^n$.

\begin{exa}\label{e.Z12}
Consider
$Z_1 = \cM \times \{0\}$ and $Z_2 = \{(x,tx) \mid x \in \cM\}$ in $K = \Hahn$
\pict[b]{15}{7cm 9cm 5.5cm 5cm}{Example~\ref{e.Z12}}(\figref). Both, $Z_1$ and $Z_2$ are $(K \times \{0\})$-riso-trivial individually (on $\cM^2$), since they are graphs of functions satisfying \eqref{eq.c1t}, but together, they are
not.\footnote{One might now think that $(Z_1, Z_2)$ is $V$-riso-trivial if and only if $Z_1 \cup Z_2$ is. But a simple counter-example consists of $Z_1 = \cM \setminus \{0\}$ and $Z_2 = \{0\}$.}
\end{exa}

By the way, Example~\ref{e.Z12} suggests that
for $V_1 = (K \times \{0\})$ and $V_2 = \{(x, tx) \in K \mid x \in K\}$, $V_1$-riso-triviality and $V_2$-riso-triviality are equivalent conditions. This raises the following question:
\begin{exe}\label{e.rtsp}
Under which condition on sub-vector spaces $V_1, V_2 \subset K^n$ is $V_1$-riso-triviality equivalent to $V_2$-riso-triviality?

Hint: If you manage to define a natural map $\res$ sending sub-vector spaces of $K^n$ to sub-vector spaces of $k^n$, then your condition should become equivalent to $\res V_1 = \res V_2$.
\end{exe}
This exercise is crucial for the theory: it means that by considering $\res V$ instead of $V$, we obtain a \emph{canonical} vector space associated to a set $Z \subset K^n$ and a ball $B \subset K^n$. (This canonical vector space is the \emp{$\operatorname{rtsp}_B(Z)$} mentioned in Footnote~\ref{f.rtsp}.) However, we do not need this in the remainder of these notes, so feel free to skip solving the exercise.

\subsection{Riso-triviality versus $C^1$-manifolds}

Using riso-triviality, we can now reformulate
Proposition~\ref{p.mnf1} in a nicer way. Let us formulate it right away in arbitrary dimension:

\begin{prop}\label{p.mnf2}
An $L$-definable set $X \subset \RR^n$ is a $d$-dimensional $C^1$-manifold if and only if, for every $\atup \in X$, $X^*$ has riso-triviality-dimension $d$ on the infinitesimal neighbourhood $B(\atup, > 0)$.
\end{prop}

The attentive reader might have noticed a difference between the  condition of Proposition~\ref{p.mnf2} (in the case $n =2, d=1$) and the one from Proposition~\ref{p.mnf1}. Couldn't it be that, for some infinitesimal neighbourhood $B$ of a point of $X$, the risometry $\alpha\colon B \to B$ sends $X^* \cap B$ to the union of several parallel lines?
For arbitrary valuative balls $B$, this can happen (as we shall see in the next example), but using nsa-methods once more, one can check that this cannot happen if $B$ is an infinitesimal neighbourhood of the form $B = B(\atup, >0)$, for some $\atup \in \RR^n$.

\subsection{An example: The cusp curve}
\label{s.cusp}

Let us consider a concrete example in some more detail: Let $X \subset \RR^2$ be the cusp curve, defined by $x^3 - y^2 = 0$, and let us determine the riso-triviality dimension of $X^*$ on various valuative balls
\pict[b]{16}{6cm 7cm 4cm 3cm}{The cusp}(\figref):

\begin{enumerate}
 \item For $\atup \in X \setminus \{0\}$ and $B := B(\atup, >0)$, we have $\rtrdim_B(X^*) = 1$, since $X$ is a $C^1$-manifold in a neighbourhood of $\atup$. The same holds for $B := B(\atup, >\lambda)$, for any $\lambda \ge 0$.
 \item For $B := B((0,0), > \lambda)$, we have $\rtrdim_B(X^*) = 0$. (This is true for arbitrary $\lambda \in \QQ$.) 
 \item For any valuative ball $B \subset \Rs^2$ disjoint from $X$, we obviously have $\rtrdim_B(X^*) = 2$.
 \item The most interesting balls are the ones which do not contain the origin, but which meet two branches of $X^*$. On such a ball, the riso-triviality-dimension is also $1$; let us analyse this in more detail.
\end{enumerate}
For simplicity, we fix one such ball, say $B := B((t^4, 0), >4)$.
Then $X^* \cap B$ is the union of the graphs $Z_+ :=  \gr g_+$,  $Z_- := \gr g_-$ of the two functions $g_{\pm}\colon x \mapsto \pm\sqrt{x^3}$ (where $x$ runs over $B_1 := B(t^4, >4)$). Both functions satisfy \eqref{eq.c1t}, so we already know that there exist risometries $\alpha_{\pm}\colon B \to B$ sending $Z_+$ and $Z_-$ to $B_1 \times \{0\}$, respectively. (In case you did not solve Exercise~\ref{exe.c1t} (3): define $\alpha_{\pm}(x,y) := (x,y-g_{\pm}(x))$.) Our goal is to glue $\alpha_+$ and $\alpha_-$ together to a single risometry $\alpha\colon B \to B$ making both graphs translation invariant. To this end, we first need to shift $\alpha_+$ and $\alpha_-$ in such a way that the images $\alpha_+(Z_+)$ and $\alpha_-(Z_-)$ have the right vertical distance: Let us modify $\alpha_+$ so that it sends $Z_+$ to $B_1 \times \{g_+(t^4)\} = B_1 \times \{t^6\}$, and analogously let $\alpha_-$ send $Z_-$ to $B_1 \times \{-t^6\}$
\pict[b]{17}{2.5cm 9cm 4cm 5cm}{Gluing risometries}(\figref).

Now we can already set $\alpha(\atup) := \alpha_+(\atup)$ if $\atup \in Z_+$ and $\alpha(\atup) := \alpha_-(\atup)$ if $\atup \in Z_-$; but what do we do in between?
 
Solution 1: Interpolate linearly. It is not difficult to check that this works wonderfully (exercise)\dots in this particular example. However, at some point deep in some proof, we need a similar kind of gluing in a more complicated situation, where it is completely unclear (at least to me) how to even define interpolation properly. Therefore, here is another approach:

Solution 2: Cut and paste: Let $B_\pm$ be the smallest balls such that the rectangles $B_1 \times B_\pm$ contain $Z_\pm$. (An easy computation yields $B_\pm = B(\pm t^6, >6)$.) We then define
$\alpha$ to be equal to $\alpha_+$ on $B_+$, equal to $\alpha_-$ on $B_-$, and equal to the identity everywhere else.
It may be a bit surprising that an $\alpha$ glued together in such a naive way is a risometry again. The key is the following lemma:

\begin{lem}
If $\alpha \colon B \to B$ is a risometry, $B' \subset B$ is a sub-ball, and $\alpha'\colon B' \to \alpha(B')$ is a risometry from $B'$ to the image of $B'$ under $\alpha$, then we obtain a new risometry by taking $\alpha$ outside of $B'$ and $\alpha'$ inside of $B'$.
\end{lem}

(The proof is a simple exercise. Note that pictures are misleading: Two points within $B'$ are always much closer to each other than to any point outside of $B'$.)

The lemma can also be applied to infinitely many disjoint sub-balls. In particular, we can apply it in our above situation, by considering each rectangle $B_1 \times B_\pm$ as a (disjoint) union of balls of the form $B((a,\pm t^6), >t^6)$.


\subsection{Better risometries}

There is yet another lesson we can draw from
Exercise~\ref{exe.c1t}. Observe that we proved that if there exists a risometry $\alpha\colon \cM^2 \to \cM^2$ sending $Z$ to $Z' := \cM \times \{0\}$, then $Z$ satisfies \eqref{eq.tf}; and if $Z$ satisfies \eqref{eq.tf}, then we constructed a risometry $\alpha'\colon \cM^2 \to \cM^2$ sending $Z$ to $Z'$. But note that the newly constructed $\alpha'$ has the additional property that it preserves the $x$-coordinate (at least if you defined $\alpha'$ in the way I had in mind, namely $\alpha'(x,y) := (x, y - g(x))$. A similar observation works in general:

\begin{lem}\label{l.aba}
If a set $Z \subset K^n$ is $d$-riso-trivial on a valuative ball $B$, then there exists a risometry $\alpha\colon B \to B$ witnessing this (i.e., such that $\alpha(Z \cap B)$ is $d$-dimensionally translation invariant) which preserves $d$ of the coordinates, i.e., such that for a projection $\pi\colon K^n \to K^d$ to a suitable subset of the coordinates, we have
$\pi \circ \alpha = \pi$.
\end{lem}

While this lemma is a bit technical, it is extremely useful. For simplicity, I will illustrate this in the example where $Z \subset K^2$ is $V$-riso-trivial on $\cM^2$
for some $1$-dimensional $V \subset K^2$ which is ``not too steep''.
(The precise condition is that $V$ is a line whose slope has non-negative valuation.)
Suppose that $\alpha\colon \cM^2 \to \cM^2$ witnesses this $V$-riso-triviality and preserves the $x$-coordinate.
Given $x \in \cM$, let us write $F_x := \{x\} \times \cM$ for the fiber above $x$. Then firstly, $\alpha$ induces a risometry $\alpha_x \colon F_x \to F_x$ for every $x$, and secondly, using translation invariance of $\alpha(Z \cap \cM^2)$, for any $x, x' \in \cM$, we can easily specify a risometry $\beta\colon F_x \to F_{x'}$ sending $\alpha(Z \cap F_x)$ to $\alpha(Z \cap F_{x'})$
\pictl{18}{3.5cm 8cm 3cm 4cm}{Lemma~\ref{l.fib} (1)}(\figref).
Putting those maps together, we obtain a risometry $\alpha_{x'}^{-1} \circ \beta \circ \alpha_x$ from $F_x$ to $F_{x'}$ sending $Z \cap F_x$ to $Z \cap F_{x'}$. Thus, all vertical fibers of $Z$ within $\cM^2$ have the same risometry type.\mpp{risometry types}\footnote{By ``$Z \subset B$ and $Z' \subset B'$ have the same risometry type'', I mean that there exists a risometry $\alpha\colon B \to B'$ with $\alpha(Z) = Z'$.} One can then verify that many properties of $Z \cap \cM^2$ are equivalent to corresponding properties of $Z \cap F_x$. As you will have guessed, this also works in higher dimension. Here are some precise statements:

\begin{lem}\label{l.fib}
If a set $Z \subset K^n$ is $V$-riso-trivial on a ball $B$, for some $d$-dimensional $V \subset K^n$, then for a suitable coordinate projection $\pi\colon K^n \to K^d$ (where suitability only depends on $V$), we have:
\begin{enumerate}
 \item The risometry type of a fiber $Z_{\atup} := Z \cap B \cap  \pi^{-1}(\atup)$, for $\atup \in \pi(B)$, does not depend on the choice of $\atup$.
 \item The risometry type of $Z \cap B$ is determined by $V$ and the risometry type of a fiber $Z_{\atup}$.
 \item Given a sub-ball $B' \subset B$ and some $d' \ge d$, the set $Z$ is $d'$-riso-trivial on $B'$ if and only if the fiber $Z_{\atup}$ is $(d'-d)$-riso-trivial on the corresponding\footnote{For simplicity, we may without loss choose $\atup \in \pi(B')$; then the ``corresponding ball'' is $B' \cap \pi^{-1}(\atup)$.} ball in the fiber
\pictr{19}{3cm 8cm 5cm 4cm}{Lemma~\ref{l.fib} (3)}(\figref).
\end{enumerate}
\end{lem}

I formulated the lemma for a ``suitable coordinate projection'' only to keep things simple. In reality, it also works for other projections, and one can easily
make precise which projections work.

\section{Back to stratifications}
\label{s.back}

\subsection{A first trial (already pretty good)}
\label{s.trial}

Let us now use riso-triviality to define stratifications. As in Definition~\ref{d.rt}, we will work with tuples of sets, so suppose we are given a tuple $\Xtup = (X_1, \dots, X_\ell)$ of $L$-definable sets $X_i \subset \RR^n$. We want to have a single stratification for the entire tuple; for this to make sense, we will stratify the ambient space $\RR^n$ in a way which is ``compatible'' with each $X_i$.

As the reader by now probably guessed, we would like to define our stratification as follows. (We call it ``shadow'', because, as we shall see later, it is only the shadow of something much stronger.)

\begin{defn}\label{d.sh}
The \empp{shadow}{the shadow} of $X_1, \dots, X_\ell$ is the partition of $\RR^n$ into sets $S_0, \dots, S_n$ defined by
\[
S_d := \{\atup \in \RR^n \mid \rtrdim_{B(\atup, >0)}(X_1^*, \dots, X_\ell^*) = d\}.
\]
\end{defn}

Note that above the definition, I only wrote ``would like to define''. Indeed, the definition still needs one tweak, to solve a problem with one of the wishes.
Let us check which of the Wishes~\ref{w} (and \ref{w2}) the definition already satisfies:

\medskip\noindent\emph{Wish (1).} We wanted $S_d$ to be smooth. As we saw in Proposition~\ref{p.mnf2}, we only obtain that it is a $C^1$-manifold.
This is fine, I am happy with that. (It might be possible to artificially refine the stratification to make the pieces $C^\infty$, but I prefer a slightly weaker but natural notion.)

\medskip\noindent\emph{Wish (2).} We want $S_d$ to be either empty or have dimension $d$.\footnote{One might feel better to first make sure that $S_d$ is definable (Wish (4')) so that dimension makes sense, but we can also think of the manifold dimension.} Indeed, one of the big theorems of \cite{iBW.can} is:
 
\begin{thm}\label{t.dim}
$S_d$ is either empty or a $d$-dimensional $C^1$-manifold.
\end{thm}

Using nsa-methods as in the proof of Proposition~\ref{p.mnf2}, it is not so hard to see that if $S_d$ is non-empty, then it must have dimension at least $d$. The only difficulty consists in proving that $S_d$ has dimension at most $d$.

It is well known that an $L$-definable set $X$ is a $C^1$-manifold almost everywhere, so for $d_0 = \dim X$, it is clear that $\dim S_{d_0-1} \le d_0-1$. However, for $d \le d_0-2$, the theorem is not clear at all. Recall that in Example~\ref{exa.crr}, we saw a set where, if we would use smooth $d$-triviality, we would obtain a one-dimensional $S_0$. The difficulty of Theorem~\ref{t.dim} consists in showing that such things do not happen with riso-triviality.

Proving Theorem~\ref{t.dim} boils down to similar questions as proving the existence of Whitney or Lipschitz stratifications. I will not give any details in these notes, but I will relate it to a result which talks purely about valued fields (Corollary~\ref{c.Trd}).

\medskip\noindent\emph{Wish (3).} We want $S_0 \cup \dots \cup S_d$ to be topologically closed. Unfortunately, this can fail. Before considering an example (and also a solution to the problem), let us consider the remaining wishes.

\medskip\noindent\emph{Wish (4').} We want $S_d$ to be $L$-definable. Indeed, another big theorem of \cite{iBW.can} is:

\begin{thm}\label{t.def}
$S_d$ is $L$-definable.
\end{thm}

In Section~\ref{s.def}, I will present some of the key techniques used in the proof.

\medskip\noindent\emph{Wish (5).} We want the stratification to be as strong as possible; more precisely, at least the origin in Example~\ref{exa.whit} should be in $S_0$. This is the case. More generally, one can verify (using nsa-methods) that our $S_0, \dots, S_n$ satisfies Whitney's regularity conditions (by \cite[Theorem~4.6.4]{iBW.can}).

\subsection{Solving the topology problem}

Now back to Wish (3):

\begin{exa}\label{e.to}
Let $X \subset \RR^3$ be the graph of the function
\[
f(x,y) := \begin{cases}
           0 & \text{if } x \le 0\\
           x\cdot y & \text{if } x > 0.
          \end{cases}
\]
(Recall that we have the order in our language $L$, so $f$ is $L$-definable.) As one can guess from the picture
\pict[b]{20}{2.5cm 7cm 3.5cm 2.5cm}{Example~\ref{e.to}}(\figref),
$S_1$ contains the $y$-axis, except maybe the origin.
On the infinitesimal neighbourhood $\cM^3$ of the origin however, $f^*$ satisfies a two-dimensional version of \eqref{eq.c1t}, so that $\rtrdim_{\cM^3}(X^*) = 2$. This means that the origin lies in $S_2$. In particular, $S_0 \cup S_1 = S_1$ is not topologically closed. 
\end{exa}

To get rid of this problem, let $(S'_0, S'_1, S'_2, S'_3)$ be the shadow of the shadow $(S_0, S_1, S_2, S_3)$.
Since $S_1 = \{(0,y,0) \in \RR^3 \mid y \ne 0\}$, we have 
$S^*_1 = \{(0,y,0) \in \Rs^3 \mid y \ne 0\}$, so the tuple of sets $(S^*_0, S^*_1, S^*_2, S^*_3)$ has riso-triviality-dimension $0$ on $\cM^3$, and hence $(0,0,0) \in S'_0$. (The rest remains unchanged: $S'_1 = S_1$, $S'_2 = S_2 \setminus \{(0,0,0)\}$, $S'_3 = S_3 = \RR^3 \setminus X$.)

One can verify that in general (i.e., if we start with an arbitrary tuple of $L$-definable sets in $\RR^n$), after sufficiently many iterations of taking the shadow, the result stabilizes, and this stabilized shadow satisfies Wish~(3). More precisely, $n$ iterations are always enough. So here is the final definition:

\vbox{%
\includegraphics[page=21,trim={0cm 12cm 0cm 4.5cm},clip,scale=0.6]{pictures}
\begin{defn}
Given a tuple $\Xtup = (X_1, \dots, X_\ell)$ of $L$-definable sets $X_i \subset \RR^n$, take its shadow, take the shadow of that, etc. The $n$-th iteration is called the \emph{riso-stratification} of $\Xtup$.
\end{defn}
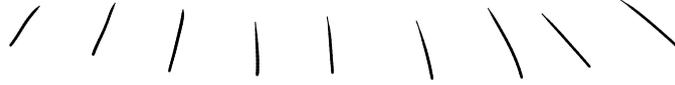
\includegraphics[page=21,trim={0cm 7.5cm 0cm 9cm},clip,scale=0.6]{pictures}
}\mpp{the riso-stratification}

We made sure that Wish~(3) becomes true. One easily verifies (exercise) that iterating does not harm Wishes~(1), (2) and (4'), and it is not very difficult to prove that more iterations can only improve Wish~(5), in the sense that if a point $\atup$ lies in $S_d$ for some iteration of taking the shadow, then in later iterations, it can only lies in $S_{d'}$ with $d' \le d$.

So now we have it: Given an $L$-definable set $X \subset \RR^n$ and a point $\atup \in X$ (or even $\atup \in \RR^n$), we have a properly working notion of the ``dimension of almost translation invariance'' of $X$ near $\atup$. We can consider this number as an invariant associated to the singularity of $X$ at $\atup$, but seen like this, it seems rather weak, given that it is just a number between $0$ and $n$. Before having a look at the proof of Theorem~\ref{t.def}, I will explain how, using the notion of riso-triviality, one can gain much more information about singularities.

Let me end this section by commenting on the relation of riso-stratifications to some more classical notions of stratifications:

As stated in Wish~(5) at the end of Subection~\ref{s.trial}, the riso-stratification of a set is in particular a Whitney stratification, which has the additional benefit of being canonically defined. Canonically defined Whitney stratifications can be obtained much more easily: in many settings, there is always a minimal one. However, in general, the minimal Whitney stratification is weaker: \mpp{riso-stratifications vs.\ Whitney stratifications}\cite[Example~5.4.6]{iBW.can} describes an algebraic set where the riso-stratification captures a singularity not seen by Whitney stratifications (and not by Verdier stratifications either).\footnote{Speaking of Verdier stratifications: I do \emph{not} believe that every riso-stratification satisfies Verdier's regularity condition.} Capturing this singularity is crucial for the example application to Poincaré series given in Section~\ref{s.p} below.

I do not know the precise relation between the riso-stratification and Lipschitz stratifications. Certainly, the riso-stratification is not, in general, a Lipschitz stratification, since the latter ones are inherently non-canonical (see Example~\ref{e.trumpet}). However, Example~\ref{e.trumpet} and also \cite[Theorem~4.2.6]{iC.aas} suggest that riso-stratifications are not so far from Lipschitz stratifications.\mpp{riso-stratifications vs.\ Lipschitz stratifications} In the other direction,
I believe that every Lipschitz stratification ``is riso-regular'', in the sense that if a point $\atup\in \RR^n$ lies in the $d$-dimensional stratum of a Lipschitz stratification of a set $X$, then $\rtrdim_{B(\atup,>0)}(X) \ge d$. However, I did not try to write down a proof of this claim. (This is closely related to \cite[Question~7.3.2]{iC.aas}.)

\section{The riso-tree}
\label{s.rt}

\subsection{Riso-triviality dimensions everywhere}

Let us go back to the more general setting of an arbitrary spherically closed valued field $K$ of equi-characteristic $0$.
First, note that the set of valuative balls in $K^n$ forms a tree
\pict[b]{22}{3cm 7cm 1.5cm 3.5cm}{Definition~\ref{d.Bn}}(\figref).
In that tree, the branches going upwards to $\infty$ correspond exactly to the points of $K^n$.

\begin{defn}\label{d.Bn}
Let \emp{$\cB^{(n)}$} be the set of all open and all closed valuative balls in $K^n$.
\end{defn}

Now fix a set $Z \subset K^n$. This yields a partition of $\cB^{(n)}$ according to the riso-triviality dimension:

\begin{defn}
The \empp{riso-tree}{the riso-tree $\Tr_0, \dots, \Tr_n$} of a set $Z \subset K^n$ is the partition of $\cB^{(n)}$ into the sets $\Tr_0, \dots, \Tr_n$ defined by
\[
\Tr_d := \{B \in \cB^{(n)} \mid \rtrdim_B Z = d\}.
\]
\end{defn}

Clearly, in the case $K = \RR^*$ and if $X$ is an $L$-definable subset of $\RR^n$, then the riso-tree of $X^*$ contains all the information about the shadow $(S_0, \dots, S_n)$ of $X$.
We will now consider a few examples showing that it contains a lot more interesting information about the singularities of $X$. In all concrete examples, we will work in $K = \Hahn$. 
By the way, I can now properly explain the term ``shadow'': 
$(S_0, \dots, S_n)$ is (only) the shadow of the riso-tree of $X^*$.

As a warm up, let us consider the riso-tree of a finite set:

\begin{exa}\label{e.finZ}
Consider $Z = \{0, t^2, t^2+t^4, 1, 1 + t^3\}$
\pict[b]{23}{2cm 10cm 4.5cm 2cm}{Example~\ref{e.finZ}}(\figref):
On any ball $B$ containing a point from $Z$, the riso-triviality dimension is $0$; these balls are exactly the infinite branches of $\cB^{(1)}$ corresponding to the points of $Z$. On all remaining balls, the riso-triviality dimension is $1$. Thus $\Tr_0$ is a tree with finitely many branching points, and one easily verifies that the depths of the branching points correspond to the valuations of differences of points of $Z$.

$\Tr_1$ is not a single tree, but it consists of many trees growing out of $\Tr_0$. Clearly, each of those trees has infinite branchings everywhere.
\end{exa}

In general, for any set $Z \subset K^n$ and any valuative balls $B' \subset B \subset K^n$, we have $\rtrdim_{B'} (Z) \ge \rtrdim_{B} (Z)$ (since a risometry witnessing $d$-riso-triviality on $B$ can be restricted to a risometry witnessing $d$-riso-triviality on $B'$.) From this, we deduce, for every $d$:
the union $\Tr_0 \cup \dots \cup \Tr_d$ is a tree;
and: for $d \ge1$, $\Tr_d$ is a union of trees (which we call the \emph{components} of $\Tr_d$), each of which grows out of 
$\Tr_0 \cup \dots \cup \Tr_{d-1}$.

By the way: It feels a bit unnatural that $\cB^{(n)}$ consists of two different kinds of balls, namely the open and the closed ones. By the following lemma, one looses no information if one defines the riso-tree to consist only of closed balls:\footnote{And indeed, this is how the riso-tree is defined in \cite{iBW.can}.}


\begin{lem}
A set $Z$ is $d$-riso-trivial on a valuative ball $B$ if and only if it is $d$-riso-trivial on every closed valuative sub-ball of $B$.
\end{lem}

However, for these notes, I found it easier to also include the open balls in the riso-tree, since often, these are the ones one is really interested in.

\subsection{Using fibers and finding Puiseux pairs}
\label{s.fib}

Here is a first example where the riso-tree contains interesting information not (directly) visible in the riso-stratification.

\begin{exa}\label{e.curi}
Let $X \subset \RR^2$ be the cusp curve (defined by $x^3 - y^2 = 0$), and let us determine the riso-tree of $Z := X^*$
\pict{24}{5cm 7cm 0cm 3cm}{Example~\ref{e.curi}}(\figref).
We already did the computations in Subsection~\ref{s.cusp}:
\begin{itemize}
 \item 
The balls with riso-triviality dimension $0$ are exactly the ones containing the origin, so $\Tr_0$ consists of a single branch (infinite both upwards and downwards).
\item
$\Tr_1$ consists of all the remaining balls $B$ which have non-empty intersection with $X^*$, but what is its tree structure? From an abstract, combinatorial point of view, we just have infinitely many trees growing out of $\Tr_0$, and each of those trees has infinite branchings everywhere. This is not a satisfactory description, so let us do better. To this end, we pause this example to make some observations.
\end{itemize}
\end{exa}

Observe that the riso-tree of a set $Z \subset K^n$ only depends on the risometry type of $Z$: 
\begin{rem}
If $\alpha\colon K^n \to K^n$ is a risometry, $\Tr_0, \dots, \Tr_n$ is the riso-tree of $Z$, and $\Tr'_0, \dots, \Tr'_n$ is the riso-tree of $Z' := \alpha(Z)$, then $\alpha$ induces a bijection from $\cB^{(n)}$ to $\cB^{(n)}$ preserving the tree structure (exercise) which sends $\Tr_d$ to $\Tr'_d$ for every $d$ (also exercise). \end{rem}
Clearly, a restricted variant of the remark holds when $\alpha$ is a risometry between two valuative balls $B, B' \subset K^n$.

Now that we know that the riso-tree only depends on the risometry type, Lemma~\ref{l.fib} suggests that if $Z$ is $d$-riso-trivial on a ball $B \subset K^n$, then the restriction to $B$ of its riso-tree should be determined by the riso-tree of a suitable fiber $Z_{\atup}$. Indeed, one can deduce the following from Lemma~\ref{l.fib}. Denote
by $\Tr_0, \dots, \Tr_n$ the riso-tree of $Z$ ``restricted to $B$'' (meaning that we only consider sub-balls of $B$).
Fix a suitable coordinate projection $\pi\colon K^n \to K^d$ and an arbitrary $\atup \in \pi(B)$ (as in Lemma~\ref{l.fib}) and set $B_{\atup} := B \cap \pi^{-1}(\atup)$ and $Z_{\atup} := Z \cap B_{\atup}$. We consider
$Z_{\atup}$ and $B_{\atup}$ as subsets of $K^{n-d}$ and denote by $\Tr'_0, \dots, \Tr'_{n-d}$ the riso-tree of $Z_{\atup}$ restricted to the ball $B_{\atup}$.

With this notation, each ball of $\Tr'_{d'}$ (for $0 \le d' \le n - d$)
corresponds exactly to ``$d$-dimensional many'' balls of $\Tr_{d + d'}$
\pict{25}{.5cm 10.5cm 5cm 3cm}{$\Tr_1 = \lin \times \Tr'_0$} (\figref).
Let me sloppily write \mpp{$\lin^d \times \Tr'_{d'}$}$\lin^d \times \Tr'_{d'}$ for the tree $\Tr_{d + d'}$ obtained from $\Tr'_{d'}$ in this way. This notation is very handy to properly describe riso-trees. Let me illustrate this by continuing our previous example.

\begin{exa}[Continuation of Example~\ref{e.curi}]\label{e.curi.c}
We can now describe $\Tr_1$ of the cusp curve properly
\pict[b]{26}{0cm 6.5cm 0cm 3.5cm}{Example~\ref{e.curi.c}}(\figref):
\begin{enumerate}
 \item 
Let us start by considering a ball of the form
$B := B(\atup, >0)$, for $\atup = (a_1, a_2) \in X \setminus \{(0,0)\}$.
Since the next bigger ball is already $\cO^2$, which contains the origin and hence lies in $\Tr_0$, 
each such ball is the root of one of the components of $\Tr_1$.
That component is of the form $\lin \times \Tr'_0$, where
$\Tr'_0$ is the riso-tree of a fiber of $X^*$ (restricted to the corresponding ball $B' := B(a_2, >0)$). Since $X^* \cap B$ is in risometry to a single straight line, the fiber is a singleton, so this $\Tr'_0$ simply consists of one infinite branch (starting at $B'$ and going upwards to infinity).
\item
Big balls far away are similar:
For every $\lambda < 0$ and every $\btup \in \Rs^2$ with $v(\btup) = \lambda$, the ball $B(\btup, > \lambda)$ is the root of a component of $\Tr_1$, and that component has the form $\lin \times \Tr'_0$, where again, $\Tr'_0$ is one infinite branch.
\item
The interesting things happen near the singularity, more precisely for balls of the form $B := B((a_1,0), > \lambda)$, when $a_1$ positive infinitesimal. We choose $\lambda := v(a_1)$, so that $B$ is the biggest ball not containing the origin, and hence is again the root of a component of $\Tr_1$.
That component is again of the form $\lin \times \Tr'_0$ where $\Tr'_0$ comes from a fiber, but this time, as we saw in Subsection~\ref{s.cusp}, the fiber consists of two points of valuative distance $\frac32\lambda$. Therefore, $\Tr'_0$ consists of two infinite branches, starting together at $B(0, > \lambda)$ and separating at depth $\frac32\lambda$. 
\end{enumerate}
This is all of $\Tr_1$.

Finally, $\Tr_2$ is boring again (but we can nevertheless describe it nicely using our approach using fibers): each of its components is of the form $\lin^2 \times \Tr'_0$, where $\Tr'_0$ is a single infinite branch.
\end{exa}

The interesting part of the above example is the $\frac32$ which one can read off from $\Tr_1$. This obviously comes from the exponents in the defining equation $x^3 - y^2$, and I am pretty sure that if one considers a curve in $\CC(\!(t^{\QQ})\!)^2$, one can recover all its Puiseux pairs from its $\Tr_1$ (though I never did the computation formally).

\subsection{Discovering hidden singularities}
\label{s.hid}

Consider now the following example:

\begin{exa}\label{e.h1}
Let $X \subset \RR^2$ be the hyperbola, defined by $xy -1 = 0$.
Since $X$ is smooth, $\Tr_0$ should clearly be empty, right?
Well, we know that $X^*$ is (at least) $1$-riso-trivial on every \emph{infinitesimal} ball contained in $\cO^2$, but on $\cO^2$ itself, it doesn't look that $1$-riso-trivial. (Exercise: Verify formally that it is not, using Exercise~\ref{e.rires}.) It turns out that $\Tr_0$ consist of a single branch coming from $-\infty$ and ending at $\cO^2 = B((0,0), \ge 0)$
\pictl{27}{5cm 8cm 4.5cm 4.5cm}{Example~\ref{e.h1}}(\figref).
\end{exa}

So does this mean that the hyperbola has a singularity? The answer is: somehow, yes. To see it, let us consider a variant of the above example:

\begin{exa}\label{e.h2}
Let $Z \subset \Rs^2$ be the hyperbola defined by $xy - t = 0$
\pictr[b]{28}{2.5cm 8cm 4.5cm 4.5cm}{Example~\ref{e.h2}}(\figref). I leave it to the reader to verify that this time, $\Tr_0$ is a branch coming from $-\infty$ and ending at $B((0,0), \ge \frac12)$. Now consider the image of $Z \cap \cO^2$ under the residue map. It is defined by the residue of our above polynomial, which is (since $\res t = 0$), simply equal to $xy$.
So here it is, the missing singularity, sitting at the origin of $\RR^2$.
\end{exa}

Note that while this particular $Z$ is not equal to $X^*$ for any $L$-definable $X \subset \RR^2$ (due to the $t$ appearing in its definition), it does arise naturally when working over $\RR$, namely when one considers the family of curves $xy - t = 0$, parameterized by $t$, and when one lets $t$ approach $0$.

Now we can also explain the $\Tr_0$ of Example~\ref{e.h1}. In general, any closed valuative ball $B = B(\atup, \ge \lambda) \subset K^n$ can be sent to $\cO^n$ by scaling and translating using the map
$\beta\colon B \to \cO^n, \xtup \mapsto (\xtup - \atup)\cdot t^{-\lambda}$; after that, we can apply the residue map. The riso-tree of a set $Z \subset K^n$ sees the singularities of the images of $Z$ of the form $\res(\beta(Z \cap B))$ for all such maps $\beta$.

\subsection{Seeing germs in a canonical way}

As a last example, consider the following trumpet:

\begin{exa}\label{e.trumpet}
Let $X \subset \RR^3$ be the trumpet defined by $y^2+z^2-x^3 = 0$
\pictl{29}{4.5cm 8.5cm 2.5cm 4cm}{Example~\ref{e.trumpet}}(\figref).
One easily verifies that its riso-stratification consists of $S_0 = \{(0,0,0)\}$, $S_1 = \emptyset$, $S_2 = X \setminus S_0$ and $S_3 = \RR^3 \setminus X$. But in the picture, there is more to see: just to the right of the origin, $X$ is almost translation invariant in the $x$-direction. The riso-trees captures this information, namely, if $a \in \RR^*$ is positive infinitesimal and $B$ is a ball around $(a,0,0)$ of suitable radius (e.g.\ $B = B((a,0,0), > v(a))$), then $\rtrdim_B(X^*) = 1$.
\end{exa}

This example is remarkable for the following reason: This kind of ``germ of $1$-triviality'' is also detected by Mostowski's Lipschitz stratifications, in the following way: For a stratification of the trumpet to satisfy Mostowski's regularity conditions, one is forced to include a stratum $S_1$ that leaves the origin to the right. However, this $S_1$ is inherently non-canonical. It could for example be the intersection of the trumpet with any plane containing the $x$-axis
\pictr{30}{4.5cm 8.5cm 3cm 4cm}{Uncanonical $S_1$}(\figref).
In contrast, the riso-tree provides a way to see the ``germ of $1$-triviality'' without introducing any non-canonicity.

\subsection{Riso-trees as stratifications}

I think of the riso-tree as being a kind of stratification of $\cB^{(n)}$. To justify this, one should prove that $\Tr_d$ is $d$-dimensional, for a suitable notion of dimension of trees. 
Here, one should think of ``horizontal'' dimension. To start with, the dimension condition on $\Tr_0$ is that it is finite at each level. Let me reformulate this a bit differently:

\begin{thm}\label{t.Tr0}
Suppose that $K$ is a spherically complete valued field of equi-characteristic $0$ and that
$\Tr_0, \dots, \Tr_n$ is the riso-tree of an $L'$-definable set
$Z \subset K^n$. Then there exists a finite set $Y_0 \subset K^n$ such that every ball $B \in \Tr_0$ contains at least one element of $Y_0$.
\end{thm}

The reader may have noted that I wrote $L'$ instead of $L$. Indeed, it will soon be handy to have this theorem for a more expressive language than $L$: in the definition of $Z$, we want to allow parameters from all of $K$, and we also want to allow using the valuation. Thus:

\begin{defn}\label{d.L'}
Let \mpp{the language $L'$}$L' := L_{\mathrm{ring}} \cup K \cup \{\cO\}$ be the language consisting of the ring language
$L_{\mathrm{ring}} = \{0, 1, +, -, \cdot\}$, one constant for each element of $K$, and a predicate $\cO$ for the valuation ring.
\end{defn}

\begin{conv}
As usual, we allow ourselves to replace $L'$ by any other language with the same expressive power. In particular, in $L'$-formulas, we will freely use the value group $\Gamma$ and the residue field $k$ as sorts\footnote{meaning that we can use variables running over $\Gamma$ and $k$}, the valuation and the residue map, the ordered abelian group language $\{0, +, - , <\}$ on $\Gamma$ and the ring language on $k$.
\end{conv}

%
Note that even though $L$ contains the order $<$ and $L'$ does not, every $L$-definable set is also $L'$-definable, since the order was not necessary anyway: $x \le y$ can also be expressed as $\exists z\colon y - x = z^2$. (While working in $\RR$ and $\RR^*$, it was convenient to have the order in $L$; now, we want to work more generally in valued fields $K$ which might have no order.)

Theorem~\ref{t.Tr0} is only about $\Tr_0$, but using that for each $d$, each component of $\Tr_d$ is of the form $\lin^d \times \Tr'_0$ for some suitable $\Tr'_0$, we deduce that $\Tr_d$ is of dimension $d$ in the following sense:

\begin{cor}[of Theorem~\ref{t.Tr0}?]\label{c.Trd}
In the setting of Theorem~\ref{t.Tr0}, we have, for every $d$, a $d$-dimensional $L'$-definable set $Y_d \subset K^n$ such that every ball $B \in \Tr_d$ contains at least one element of $Y_d$.
\end{cor}

Actually, this does not follow from Theorem~\ref{t.Tr0} in the way it is formulated above: While Theorem~\ref{t.Tr0} allows us to find a finite set $Y_0$ for the $\Tr'_0$ of each fiber of each component of $\Tr_d$ (and then define $Y_d$ to be the union of all those sets), we need to be able to do this in a uniform way, so that $Y_d$ becomes $L'$-definable. More precisely, we need (and have) the following:

\begin{add}[to Theorem~\ref{t.Tr0}]\label{a.Tr0}
Suppose that we have an $L'$-definable family of sets $Z_{\ctup} \subset K^n$, parameterized by $\ctup \in K^m$. (This means that there is a single $L'$-formula $\phi(\xtup, \ytup)$ such that $Z_{\ctup} = \{\atup \in K^n \mid K \models \phi(\atup, \ctup)\}$ for every $\ctup$.) Then there exists an $L'$-definable family of sets $Y_{0,\ctup}$ such that for every $\ctup$, the set $Y_0 := Y_{0,\ctup}$ satisfies the condition from Theorem~\ref{t.Tr0} applied to $Z := Z_{\ctup}$: $Y_{0,\ctup}$ is finite, and every ball of the $\Tr_0$ of $Z_{\ctup}$ contains at least one element of 
$Y_{0,\ctup}$.
\end{add}

While Corollary~\ref{c.Trd} does follow from this improved Theorem~\ref{t.Tr0}, it is still a bit dishonest to call it a corollary, since in reality, the proof goes the other way round: To prove Theorem~\ref{t.Tr0}, one does a downward induction, first constructing $Y_n$, then $Y_{n-1}$, etc., until $Y_0$, which is the most difficult to obtain.

\subsection{Deducing Theorem~\ref{t.dim}}

Recall that concerning the riso-stratification $S_0, \dots, S_n$ of an $L$-definable set $X \subset \RR^n$, the claim of Theorem~\ref{t.dim} was that $\dim S_d = d$.
Corollary~\ref{c.Trd} is the corresponding statement for riso-trees, and from that corollary, it is not too hard to deduce Theorem~\ref{t.dim}: For $d = 0$, this is clear: By the corollary, every point of $S_0$ has at least one point of $Y_0$ in its infinitesimal neighbourhood, so since $Y_0$ is finite, $S_0$ is finite, too.

For $d \ge 1$, one needs an argument to relate the dimension in $\Rsn$ (of the set $Y_d$) to dimension in $\RR^n$ (of the set $S_d$); this is not entirely straight forward, since $Y_d$ is not $L$-definable, but it is not too hard, either
(using an idea explained in Subsection~\ref{s.dSd}). The real work, however, lies in the proof of Corollary~\ref{c.Trd}.

\section{Definability}
\label{s.def}

\subsection{Making risometries accessible to model theory}

We would like to apply model theoretic methods to the riso-trees. To this end, we need them to be definable, but what do we mean by this? 
We continue working in a spherically closed valued field $K$ of equi-characteristic $0$, and we continue to use our language $L'$ from Definition~\ref{d.L'}.

\begin{tenthm}\label{tt}
Given an $L'$-definable set $Z \subset K^n$, a valuative ball $B \subset K^n$, and a natural number $d \le n$, we can express by an $L'$-formula whether $Z$ is $d$-riso-trivial on $B$.
\end{tenthm}

While intuitively, the reader may have guessed what this statement is supposed to express, formally it makes no sense: Since $Z$ and $B$ and $d$ are fixed, the $L'$-formula does not depend on anything, so it has no free variables:
Either $Z$ is $d$-riso-trivial on $B$; then take the formula to be (always) true; or it is not; then take the formula to be (always) false. What we really want to express is that Tentative Theorem~\ref{tt} works uniformly if we let some parameter vary, in a similar style as in Addendum~\ref{a.Tr0}. Here is a way to make this precise:

\addtocounter{defn}{-1}
\begin{thm}\label{t.rtdef}
Suppose that $Z_{\ctup} \subset K^n$ is an $L'$-definable family of sets and that $B_{\ctup} \subset K^n$ is an $L'$-definable family of valuative balls (both parameterized by a parameter $\ctup$). Then, for every fixed natural number $d \le n$, there exists an $L'$-formula $\psi_d(\ytup)$ which expresses that $Z_{\ctup}$ is $d$-riso-trivial on $B_{\ctup}$, i.e., for every $\ctup$, we have $\rtrdim_{B_{\ctup}}(Z_{\ctup}) \ge d$ if and only if $K \models \psi_d(\ctup)$.
\end{thm}

Let me give a baby example of an application of this theorem, to illustrate that this formulation is really what one wants.
%
In our trumpet example (Example~\ref{e.trumpet}), we might want to express that if we fix an $x$-coordinate, then there exists a valuative radius such that the fiber $X^*_x$ of the trumpet at $x$ (considered as a subset of $\Rs^2$) is $1$-riso-trivial on every valuative ball of that radius:
\[
\forall x \in \RR^*\colon \exists \lambda\in \Gamma\colon
\forall (y,z) \in \Rs^2\colon
\underbrace{\text{$X^*_x$ is $1$-riso-trivial on $B((y,z), > \lambda)$}}_{(\star)}.
\]
Formally, the part denoted by $(\star)$ is not an $L'$-formula. However, in that part, both $X^*_x$ and $B((x,y,z), > \lambda)$ are given by a formula with parameters $(x,y,z,\lambda)$, so
Theorem~\ref{t.rtdef} states that $(\star)$ is equivalent to some $L'$-formula $\psi(x,y,z,\lambda)$. In other words, the theorem allows us to ``without loss'' consider statements like $(\star)$ as $L'$-formulas (and use them within other $L'$-formulas).

Theorem~\ref{t.rtdef} is proved in a common induction with another result, which is interesting in its own right, namely that having the same risometry type can be expressed in $L'$. This time, I give a correct formulation right away:

\begin{thm}\label{t.ridef}
Suppose that $Z_{\ctup}, Z'_{\ctup} \subset K^n$ are $L'$-definable families of sets and $B_{\ctup}, B'_{\ctup} \subset K^n$ are $L'$-definable families of valuative balls. Then there exists a formula $\psi(\ctup)$ which expresses that there exists a risometry $\alpha_{\ctup}\colon B_{\ctup} \to B'_{\ctup}$ with $\alpha_{\ctup}(Z_{\ctup}) = Z'_{\ctup}$.
\end{thm}

Note that while many objects in that theorem are $L'$-definable, nothing is said about definability of $\alpha_{\ctup}$. (One can construct pathological examples where there exists only a non-definable risometry sending $Z_{\ctup}$ to $Z'_{\ctup}$; in that case, $\psi(\ctup)$ will still be true.)

I will not go into details concerning the proofs of Theorems~\ref{t.rtdef} and \ref{t.ridef}, but let me at least
try to give a vague intuition about how it can be possible to express the existence of arbitrary risometries using a first order formula. Let me start as follows:

\begin{proof}[Sketch of proof of Theorem~\ref{t.rtdef} in a very simple case]
We already saw such a proof in Exercise~\ref{exe.c1t}: Being in risometry to a single horizontal line is equivalent to Condition~\eqref{eq.c1t2}, which is an $L'$-formula. It is not too difficult to generalize this a bit, replacing the single straight line by finitely many parallel lines. Moreover, similar arguments also work in higher dimension.
\end{proof}

\begin{proof}[Sketch of proof of Theorem~\ref{t.ridef} in a very simple case]
There exists a risometry between two finite sets $Z$ and $Z'$ if and only if there exist enumerations $Z = \{a_1, \dots, a_n\}$ and $Z' = \{a'_1, \dots, a'_n\}$ such that $\rv(a_i - a_j) = \rv(a'_i - a'_j)$ for every $1 \le i < j \le n$ (exercise). Such a risometry can then always be extended to ambient balls $B$, $B'$ of the same radius (also exercise).
Using this observation, it is not difficult to prove 
Theorem~\ref{t.ridef} under the assumption that the cardinalities of both, $Z_{\ctup}$ and $Z'_{\ctup}$, are bounded by some $N \in \NN$ which is independent of $\ctup$ (third exercise).
\end{proof}

So far, the reader might complain that the above two simple cases are far too special to make the theorems plausible. However, recall that we saw, in Subsection~\ref{s.fib}, that the riso-tree of a set $Z$ can be described entirely in terms of the $\Tr'_0$-parts of riso-trees of suitable fibers. By Theorem~\ref{t.Tr0}, we know that each of those $\Tr'_0$ is controlled by a finite set $Y_0$. By exploiting this cleverly, one can reduce the proof of Theorem~\ref{t.ridef} to the case of finite sets and the proof of Theorems~\ref{t.rtdef} to the case of finite fibers, i.e., to exactly the simple cases described above.
(Note however that this is just a glimpse of the proof; I am hiding a lot of work under the carpet.)

\subsection{An example application: something like higher dimensional Puiseux pairs}

For surfaces like the trumpet (Example~\ref{e.trumpet}), a natural question is: How does the valuative radius of a fiber $X^*_x$ of the trumpet depend on the valuation of the $x$-coordinate?
Let us denote the function sending $v(x)$ to that radius by $\rho$
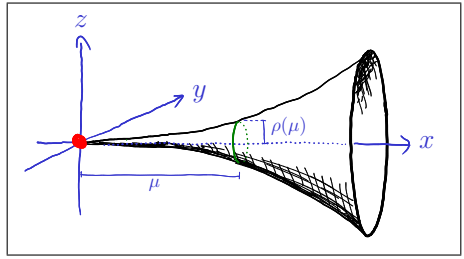
\pict{31}{4.5cm 8.5cm 3cm 4cm}{something like Puiseux pairs}(\figref).
That valuative radius is the largest $\lambda \in \Gamma$ such that $X^*_x$ has riso-triviality-dimension $0$ on $B((0,0), \ge \lambda)$, so using Theorem~\ref{t.rtdef}, we obtain that $\rho$ is $L'$-definable. (In the case of our trumpet, we have $\rho(\mu) = \frac32\mu$, but we can define a similar $\rho$ for arbitrary $L$-definable sets $X^* \subset \Rsn = (\Hahn)^n$; also, we can let $\rho$ depend on the valuations of several coordinates.)

Using general results from model theory of valued fields, one obtains that any $L'$-definable function $\rho$ from the value group of $\Rsn$ to itself is piece-wise linear. (First, one uses Denef-Pas quantifier elimination to show that $\rho$ can be defined using only the ordered abelian group language on the value group $\Gamma = \QQ$; then one concludes using quantifier elimination in divisible ordered abelian groups.)

Thus, we can for example take the slopes of our $\rho$ (which are rational numbers) and attach them as invariants to the singularity of $X$ at $0$. By this and similar methods, one can define many invariants which feel a bit like generalized Puiseux pairs, or maybe also like slopes of certain Newton polytopes. It would be interesting to relate those numbers to other invariants of singularities.

%
%

\subsection{Deducing definability of the $S_d$}
\label{s.dSd}

To end this section, let me explain how 
Theorem~\ref{t.rtdef} implies that the sets $S_d$ from our riso-stratifications are $L$-definable (Theorem~\ref{t.def}). This is not very difficult, but it uses a nice little trick.

Recall that we start with an $L$-definable set $X \subset \RR^n$ and that $S_d$ is defined as
\[
S_d := \{\atup \in \RR^n \mid \underbrace{\rtrdim_{B(\atup, >0)}(X^*) = d}_{(\star)}\}.
\]
By Theorem~\ref{t.rtdef}, we know that the set $S^*_d$ of all $\atup \in \Rsn$ satisfying $(\star)$ is definable, and $S_d$ is just the set definable by the same formula in $\RR$, isn't it? Well, no: To define $S^*_d$, we need the valuation, so the formula defining it makes no sense in $\RR$.\footnote{Equipping $\RR$ with the trivial valuation does not yield anything meaningful.} So instead of considering $\RR$ as a sub-field of $\RR^*$, we now consider it as its residue field. From this perspective, we can write $S_d$ as follows:
\[
S_d = \{\atup \in \RR^n \mid \underbrace{\rtrdim_{\res^{-1}(\atup)}(X^*) = d}_{(\star)}\}.
\]
Now $S_d$ itself is $L'$-definable, when considered as a set living in the residue field sort. While this might not yet look better, we can again use a general result from model theory of valued fields (again Denef-Pas quantifier elimination), to
obtain that every $L'$-definable subset of a Cartesian power of the residue field is already definable by an $L$-formula in the residue field. Applying this to $S_d$ yields exactly what we want.
%

\subsection{Generaliations}
\label{s.gen}

In these notes, I decided to present everything for the real numbers and using the language of ordered rings. However, everything also works in various other contexts:

\begin{itemize}
 \item Instead of taking for $L$ the language of ordered rings, one can take the sub-analytic language. This is an expansion of $L$ by function symbols for all analytic functions converging on a neighbourhood of their domain.
 In this way, one can apply the entire theory to analytic sets $X \subset \RR^n$ (and even to globally sub-analytic ones, which are exactly the $L$-definable sets); of course, the sets $S_d$ will then also only be globally sub-analytic.
 
 To see that the theory works with this expanded $L$, one 
 firstly needs to make sure that one can still find an elementary extension of $\RR$ which is spherically complete; this is indeed the case. And secondly, one needs that in the corresponding expansion of $L'$, one still has something like Denef-Pas quantifier elimination. This, too, is the case.

 (In contrast, if we expand $L$ by a function symbol for the exponential function, then -- while the model theory of $\RR$ in that language is still well-behaved -- neither of the two above requirements would be satisfied anymore.)
 \item
 Instead of working in $\RR$, we can work in $\CC$. Of course, $L$ should then be the language of rings, and not the language of ordered rings. Note that it makes a difference whether we apply the theory in $\CC$ or whether we identify $\CC$ with $\RR^2$ and just apply the theory in the structure $\RR$ as we did before: With the latter approach, to an $L$-definable set $X \subset \CC^n = \RR^{2n}$ we would associate a riso-stratification $S_0, \dots, S_{2n}$, where $S_d$ has real dimension $d$, whereas when we work directly in the structure $\CC$, the riso-stratification only consists of $S_0, \dots, S_n$, where $S_d$ has complex dimension $d$.
 
 In some places, topology plays a role, namely to make sense of small neighbourhoods and when we require that $S_0 \cup \dots \cup S_d$ is topologically closed. We use the analytic topology on $\CC$. However, for $L$-definable sets, being closed in the analytic topology is equivalent to being closed in the Zariski topology, so for the end result, the topology does not matter.
 \item
 Model theory does not distinguish between elementary equivalent structures, so instead of $\RR$ or $\CC$, we can also use arbitrary real closed fields or algebraically closed fields of characteristic zero.
 
 If we work in a real closed field $k$, we use the interval topology. While such a $k$ can be totally disconnected, this is not noticed by our theory, because we only use definable sets: Exercise: Show that any real closed field $k$ is \emph{definably} connected, i.e., that it cannot be written as the disjoint union of two non-empty $L$-\emph{definable} closed sets (for $L = L_{\mathrm{oring}} \cup k$).
 
 If we work in an algebraically closed field $k$ of characteristic $0$, the easiest way to get the theory working is by inventing an analytic topology on $k$. One can obtain such a topology by fixing a real closed sub-field $k_0 \subset k$ satisfying $[k:k_0] = 2$.
 This feels a bit awkward, since $k_0$ is nowhere near unique. However, the final results do not depend on the chosen topology, since again, for definable sets, being closed in that analytic topology is equivalent to being closed in the Zariski topology.
 \item
 Finally, people have also been interested in constructing Whitney stratifications over the $p$-adic numbers $\QQ_p$ and over fields of Laurent series like $k_0(\!(t)\!)$, where $k_0$ is an arbitrary field of characteristic $0$ \cite{CCL.cones,For.motCones}. Our theory can also be applied to those cases. Note that things start getting a bit messy here:
 We again pick an elementary extension $\QQ^*_p$ or $k_0(\!(t)\!)^*$, and the valuation we use on that elementary extension is \emph{not} the one coming from $\QQ_p$ or $k_0(\!(t)\!)$. Indeed, in terms of our valuation, 
 the infinitesimal neighbourhood of $0$ in $\QQ^*_p$ consists 
 of elements of positive valuation; but from the point of view of the $p$-adic valuation, infinitesimal means having infinite valuation (i.e., bigger than $n$ for every $n \in \ZZ$).
\end{itemize}

\section{Application to Poincaré series}
\label{s.p}

\subsection{Poincaré series in $\QQ_p$}

Let me finish these notes by sketching an application of riso-stratifications to Poincaré series.

I start by recalling what a Poincaré series is. There are several variants of that notion. The one I will present is not the most common one, but it is simpler to explain in these notes.

Let us first work in $\QQ_p$ for some time, and fix a set $Z \subset \QQ_p^n$. Note that $\ZZ_p^n$ contains only finitely many closed valuative balls of a fixed valuative radius $\lambda \in \NN$; denote by $N_\lambda$ the number of such balls which have non-empty intersection with $Z$ and use those numbers to form a formal power series $\sum_{\lambda \ge 0} N_\lambda T^\lambda$; this is the Poincaré series of $Z$.
\pict[b]{32}{3.5cm 11cm 5cm 1.5cm}{Definition~\ref{d.PZB}}(\figref).

It will be handy to also have a variant of this where instead of counting sub-balls of $\ZZ_p^n$, we count sub-balls of any fixed valuative ball $B \subset \QQ_p^n$:
\begin{defn}\label{d.PZB}
Given a set $Z \subset \QQ_p^n$ and a valuative ball $B \subset \QQ_p^n$, we define the \empp{Poincaré series}{the Poincaré series $P_{Z,B}(T)$} of $Z$ on $B$ to be the formal series
\[
P_{Z,B}(T) := \sum_{\lambda \in \ZZ} N_\lambda T^\lambda,
\]
where $N_\lambda$ is the number of closed valuative balls $B' \subset B$ of radius $\lambda$ satisfying $B' \cap Z \ne \emptyset$.
\end{defn}

Now let me explain why riso-triviality is very useful to determine Poincaré series. We allow ourselves temporarily to use riso-triviality in $\QQ_p$ even though $\QQ_p$ is not of equi-characteristic $0$. (But be reassured: we will not use any serious theorem.) Firstly, note that risometries preserve Poincaré series:

\begin{lem}\label{l.rp}
If $Z \subset \QQ_p^n$ is an arbitrary set and
$\alpha\colon B \to B'$ is a risometry, for some valuative balls $B, B' \subset \QQ_p^n$, then $P_{Z,B}(T) = P_{\alpha(Z), B'}(T)$.
\end{lem}

This follows easily by noting that $\alpha$ induces a bijection on the set of all balls of a fixed radius. (This only needs that $\alpha$ is an isometry.)

The next thing to note is that if a set is translation invariant, then its Poincaré series can easily be computed from the Poincaré series of a fiber. Indeed, suppose for example that $Z = \ZZ_p \times Z'$, for an arbitrary set $Z' \subset \ZZ_p$. Then every valuative ball $B(b, \ge \lambda) \subset \ZZ_p$ which has non-empty intersection with $Z'$ gives rise to exactly $p^\lambda$ valuative balls $B((a,b), \ge \lambda) \subset \ZZ_p^2$ which have non-empty intersection with $Z$
\pict{33}{3cm 8cm 4cm 4cm}{$P_{Z',\ZZ_p}(T)$ vs.\ $P_{Z,\ZZ^2_p}(T)$}(\figref).
Thus, if we denote by
\[
P_{Z, \ZZ_p^2}(T) = \sum_\lambda N_\lambda T^\lambda
\qquad \text{and} \qquad
P_{Z', \ZZ_p}(T) = \sum_\lambda N'_\lambda T^\lambda
\]
the Poincaré series of $Z$ and the of the fiber $Z'$, respectively, then we obtain
$N_\lambda = p^\lambda N'_\lambda$, and hence
\[
P_{Z, \ZZ_p^2}(T) = \sum_\lambda p^\lambda N'_\lambda T^\lambda
= \sum_\lambda N'_\lambda \cdot (pT)^\lambda = P_{Z', \ZZ_p}(pT).
\]
One obtains a similar formula for sets which are translation invariant in more directions (namely, $P_{Z, \ZZ_p}(T) = P_{Z', \ZZ_p}(p^dT)$), and combining this with Lemma~\ref{l.rp} yields that the Poincaré series of a $d$-riso-trivial set can be computed from the Poincaré series of a fiber:

\begin{prop}\label{p.rtQp}
If $Z\subset \QQ_p^n$ is $d$-riso-trivial on a valuative ball $B \subset\QQ_p^n$, then for a suitable coordinate projection $\pi\colon \QQ_p^n \to \QQ_p^d$ and any $\atup \in \pi(B)$, we have
$P_{Z, B}(T) = P_{Z_{\atup}, B_{\atup}}(p^dT)$, where $Z_{\atup}$ and
$B_{\atup}$ denote the $\pi$-fibers over $\atup$ of $Z$ and $B$,
respectively, considered as subsets of $\QQ_p^{n-d}$. 
\end{prop}

It should be clear now that our theory of riso-trees is useful to compute Poincaré series, except that so far, we developed the theory of riso-trees only in equi-characteristic $0$. So
since so far, the mountain did not go to Mahomet, let Mahomet go to the mountain: let me sketch how one can make sense of Poincaré series in $\Cbb{t}$.

\subsection{Motivic Poincaré series in $\Cbb{t}$}

The problem is that counting balls in the classical way makes no sense anymore, since there are almost always just uncountably many.
The solution is to use a form of ``motivic counting'' provided by motivic integration. (The idea of motivic Poincaré series goes back to \cite{DL.def}.) The idea is that instead of counting using only natural numbers, we introduce formal symbols which stand for different infinite numbers. For example, we introduce a formal symbol $[\CC]$ which stands for ``the number of elements of $\CC$''. Then we can say things like: The valuation ring
$\CC[[t]]$ contains $[\CC]^2$ many closed valuative balls of valuative radius $2$. This sounds fancy, but it is not even that difficult to make formal: Our ``numbers'' are simply elements of the Grothendieck ring $\Gro$\mpp{the Grothendieck ring of varieties $\Gro$} of varieties over $\CC$. This ring $\Gro$ is defined to contain an element $[X] \in \Gro$ for every algebraic set $X \subset \CC^n$. For our purposes, this $[X]$ stands for the ``number of elements'' of $X$.

Note that more generally, we can associate an element $[X] \in \Gro$ to every $L$-definable set $X \subset \CC^n$; here, we take $L := L_{\mathrm{ring}} \cup \CC = \{0, 1, +, -, \cdot\} \cup \CC$ to be the language of rings with constants for the elements of $\CC$ added. Indeed, by quantifier elimination in algebraically closed fields, every $L$-definable set is just a finite Boolean combination of algebraic sets, and if, for example, we have $X = X_1 \setminus X_2$ for algebraic sets $X_1, X_2 \subset \CC^n$, then $X_1 \cap X_2$ is also algebraic, and
we naturally have $[X] = [X_1] - [X_1 \cap X_2] \in \Gro$.

Using elements of $\Gro$ as numbers, it will certainly not be possible to assign a Poincaré series to an arbitrary subset of $\Cbb{t}^n$, but using Denef-Pas quantifier elimination, one
can make sense of the Poincaré series associated to any $L'$-definable set (where $L'$ is still the language from Definition~\ref{d.L'}). Here is a not entirely precise definition of that series:

\begin{propdef}
Suppose that we are given an $L'$-definable set $Z \subset \Cbb{t}^n$, a valuative ball $B \subset \Cbb{t}^n$, and a $\lambda \in \ZZ$. Consider the set of all valuative balls $B' \subset B$ of valuative radius $\lambda$ which have non-empty intersection with $Z$. This set can be parameterized by an $L$-definable set $Y_\lambda \subset \CC^{n'}$. The \empp{motivic Poincaré series}{the motivic Poincaré series $P\mot_{Z,B}(T)$} of $Z$ on $B$ is then defined to be the formal series
\[
P\mot_{Z,B}(T) := \sum_{\lambda} [Y_\lambda] \cdot T^\lambda,
\]
with coefficients $[Y_\lambda] \in \Gro$.
\end{propdef}

One is particularly interested in the following motivic Poincaré series:

\begin{defn}
Suppose we are given an algebraic set $X = X(\CC) \subset \CC^n$; let me write $X(\Cbb{t})$ for the subset of $\Cbb{t}^n$ defined by the same polynomials. Then the \empp{local motivic Poincaré series}{the local motivic Poincaré series $P\motloc_{X,\atup}(T)$} of $X$ at a point $\atup \in \CC^n$ is defined to be the motivic Poincaré series of $X$ on the infinitesimal neighbourhood of $\atup$ within $\Cbb{t}^n$:
\[
P\motloc_{X,\atup}(T) := P\mot_{X(\Cbb{t}), B(\atup, \ge1)}(T)
\]
\end{defn}
This Poincaré series $P\motloc_{X,\atup}(T)$ contains interesting (but still somewhat mysterious) information about the singularity of $X$ at $\atup$. We shall now see what riso-stratifications tell us about them.

\subsection{The connection to riso-stratifications}

One can verify that Proposition~\ref{p.rtQp} also holds in this motivic setting, if we repace the $p$ in the formula by the ``number of elements of the residue field'', i.e., by $[\CC]$:

\begin{prop}\label{p.rtCt}
If $Z \subset \Cbb{t}^n$ is an $L'$-definable set which is $d$-riso-trivial on a valuative ball $B \subset\Cbb{t}^n$, then for a suitable coordinate projection $\pi\colon \Cbb{t}^n \to \Cbb{t}^d$ and any $\atup \in \pi(B)$, we have
$P_{Z, B}(T) = P_{Z_{\atup}, B_{\atup}}([\CC]^d\cdot T)$, where $Z_{\atup}$ and
$B_{\atup}$ denote the $\pi$-fibers over $\atup$ of $Z$ and $B$,
respectively, considered as subsets of $\Cbb{t}^{n-d}$. 
\end{prop}

One might believe that once one is familiar with the motivic formalism, the proof of this proposition is just a rather straight forward adaptation of the proof of Proposition~\ref{p.rtQp}. However, this is only partially true: The problem consists in proving Lemma~\ref{l.rp} in this setting, i.e., that a risometry $\alpha$ preserves the motivic coefficients $[Y_\lambda]$ of $P_{Z, B}(T)$. This would be easy if $\alpha$ were assumed to be definable. Proving the lemma without this assumption requires analysing the $L'$-definable set $Z$ in detail using its riso-tree (and describing the tree in terms of $\Tr'_0$-parts of suitable fibers).

Once Proposition~\ref{p.rtCt} is proved, we can combine it with riso-stratifications. The result is that local motivic Poincaré series can be computed by looking just at a fiber transversal to the riso-stratification:

\begin{thm}\label{t.pc}
Suppose that $X \subset \CC^n$ is an algebraic set, let $S_0, \dots, S_n$ be its riso-stratification, and pick any point $\atup \in S_d$, for some $d$.
Then, for a suitable coordinate projection $\pi\colon \CC^n \to \CC^d$, we have
\[
P\motloc_{X,\atup}(T) = P\motloc_{X_{\atup},\atup}([\CC]^d\cdot T),
\]
where $X_{\atup} = X \cap \pi^{-1}(\pi(\atup))$ is the $\pi$-fiber of $X$ containing $\atup$.
\end{thm}


Proving this result needs one more ingredient:
Poincaré series are defined using $\Cbb{t}$, while the riso-stratification is defined using (e.g.) $\Cbb{t^{\QQ}}$. Thus, we need to verify that $d$-triviality over $\Cbb{t^{\QQ}}$ implies $d$-triviality over the sub-field $\Cbb{t}$. This is yet another proof where we need to analyse the involved definable sets using their riso-trees.

Note that Theorem~\ref{t.pc} would not be true if, instead of the riso-stratification, one would use a Whitney stratification $S_0, \dots, S_n$; an example is given in \cite{iBW.can} (namely Example~5.4.6).

\bibliographystyle{amsplain}
\bibliography{references}

\providecommand{\bysame}{\leavevmode\hbox to3em{\hrulefill}\thinspace}
\providecommand{\MR}{\relax\ifhmode\unskip\space\fi MR }
\providecommand{\MRhref}[2]{%
  \href{http://www.ams.org/mathscinet-getitem?mr=#1}{#2}
}
\providecommand{\href}[2]{#2}
\begin{thebibliography}{10}

\bibitem{iBW.can}
David Bradley-Williams and Immanuel Halupczok, \emph{Riso-stratifications and a
  tree invariant}, 2022, Preprint.

\bibitem{CCL.cones}
Raf Cluckers, Georges Comte, and Fran{\c{c}}ois Loeser, \emph{Local metric
  properties and regular stratifications of {$p$}-adic definable sets},
  Comment. Math. Helv. \textbf{87} (2012), no.~4, 963--1009. \MR{2984578}

\bibitem{iC.aas}
Pablo Cubides~Kovacsics and Immanuel Halupczok, \emph{Arc-wise analytic
  t-stratifications}, 2021, Preprint.

\bibitem{DL.def}
Jan Denef and Fran{\c{c}}ois Loeser, \emph{Definable sets, motives and
  {$p$}-adic integrals}, J. Amer. Math. Soc. \textbf{14} (2001), no.~2,
  429--469 (electronic). \MR{MR1815218 (2002k:14033)}

\bibitem{For.motCones}
Arthur Forey, \emph{Motivic local density}, Math. Z. \textbf{287} (2017),
  no.~1-2, 361--403. \MR{3694680}

\bibitem{i.whit}
Immanuel Halupczok, \emph{Non-{A}rchimedean {W}hitney stratifications}, Proc.
  Lond. Math. Soc. (3) \textbf{109} (2014), no.~5, 1304--1362. \MR{3283619}

\bibitem{iY.sts}
Immanuel Halupczok and Yimu Yin, \emph{Lipschitz stratifications in
  power-bounded {$o$}-minimal fields}, J. Eur. Math. Soc. (JEMS) \textbf{20}
  (2018), no.~11, 2717--2767. \MR{3861807}

\bibitem{Kap.maxValFlds}
Irving Kaplansky, \emph{Maximal fields with valuations}, Duke Math. J.
  \textbf{9} (1942), 303--321. \MR{6161}

\bibitem{Mos.biLip}
Tadeusz Mostowski, \emph{Lipschitz equisingularity}, Dissertationes Math.
  (Rozprawy Mat.) \textbf{243} (1985), 46. \MR{808226 (87e:32008)}

\bibitem{Whi.strat}
Hassler Whitney, \emph{Tangents to an analytic variety}, Ann. of Math. (2)
  \textbf{81} (1965), 496--549. \MR{0192520 (33 \#745)}

\end{thebibliography}

\end{document}